\documentclass[10pt,twocolumn]{IEEEtran}
\makeatletter
\def\ps@headings{%
\def\@oddhead{\mbox{}\scriptsize\rightmark \hfil \thepage}%
\def\@evenhead{\scriptsize\thepage \hfil \leftmark\mbox{}}%
\def\@oddfoot{}%
\def\@evenfoot{}}
\makeatother
\usepackage{graphicx,amsmath,amssymb,amsfonts,bm,setspace,amsthm}
\usepackage[noadjust]{cite}
\usepackage{braket} 		

\usepackage{booktabs,subfigure}
\newcommand{\otoprule}{\midrule[\heavyrulewidth]}

\allowdisplaybreaks

\theoremstyle{definition}
\newtheorem{assu}{Assumption}
\newtheorem{cor}{Corollary}
\newtheorem{defn}{Definition}
\newtheorem{lem}{Lemma}

\newtheorem{thm}{Theorem}

\newcommand{\expect}[1]{\mathbb{E}\left[#1\right]}

\newcommand{\Z}{\mathbb{Z}}

\newcommand{\argmax}{\operatorname{argmax}}

\newcommand{\abs}[1]{\left|#1\right|}

\newcommand{\amax}{A_{\text{max}}}

\newcommand{\dmax}{d_{\text{max}}}
\newcommand{\Dmax}{D_{\text{max}}}
\newcommand{\Qmax}{Q_{\text{max}}}
\newcommand{\Zmax}{Z_{\text{max}}}
\newcommand{\mumax}{\mu_{\text{max}}}

\newcommand{\ORA}{\mathsf{ORA}}
\newcommand{\UORA}{\mathsf{UORA}}

\newcommand{\Qnc}{Q_{n}^{(c)}}
\newcommand{\qnc}{q_{n}^{(c)}}
\newcommand{\Dnc}{D_{n}^{(c)}}
\newcommand{\dnc}{d_{n}^{(c)}}
\newcommand{\dncstar}{d_{n}^{(c)*}}
\newcommand{\varphinc}{\varphi_{n}^{(c)}}
\newcommand{\varphincstar}{\varphi_{n}^{(c)*}}
\newcommand{\tdnc}{\widetilde{d}_{n}^{(c)}}
\newcommand{\nuc}{\nu_{c}}

\newcommand{\numax}{\nu_{\text{max}}}
\newcommand{\deltamax}{\delta_{\text{max}}}
\newcommand{\tnu}{\widetilde{\nu}}
\newcommand{\tmu}{\widetilde{\mu}}
\newcommand{\Anc}{A_{n}^{(c)}}

\newcommand{\lambdanc}{\lambda_{n}^{(c)}}
\newcommand{\Zc}{Z_{c}}
\newcommand{\blambda}{\bm{\lambda}}

\newcommand{\bigparen}[1]{\big(#1\big)}
\newcommand{\mumaxin}{\mu_{\text{max}}^{\text{in}}}
\newcommand{\mumaxout}{\mu_{\text{max}}^{\text{out}}}
\newcommand{\onecR}{1_{c}^{\text{R}}}
\newcommand{\onecL}{1_{c}^{\text{L}}}

\newcommand{\cN}{\mathcal{N}}
\newcommand{\cL}{\mathcal{L}}
\newcommand{\cC}{\mathcal{C}}

\newcommand*{\xfill}[1][0pt]{%
   Leaders
    \hbox to 1pt{\hss
      \raisebox{#1}{\rule{1.2pt}{0.4pt}}%
      \hss}\hfill}

\begin{document}

\title{Receiver-Based Flow Control for Networks in Overload}

\author{\IEEEauthorblockN{Chih-ping Li and Eytan Modiano} \\
\IEEEauthorblockA{Laboratory for Information and Decision Systems\\
Massachusetts Institute of Technology
}
\thanks{This work was supported by  DTRA grants HDTRA1-07-1-0004 and HDTRA-09-1-005, ARO Muri grant number W911NF-08-1-0238, and NSF grant CNS-1116209.}
}

\markboth{}{}
\maketitle

\begin{abstract}
We consider utility maximization in networks where the sources do not employ flow control and may consequently overload the network. In the absence of flow control at the sources, some packets will inevitably have to be dropped when the network is in overload.  To that end, we first develop a distributed, threshold-based packet dropping policy that maximizes the weighted sum throughput. Next, we consider utility maximization and develop a receiver-based flow control scheme that, when combined with threshold-based packet dropping, achieves the optimal utility. The flow control scheme creates virtual queues at the receivers as a \emph{push-back} mechanism to optimize the amount of data delivered to the destinations via back-pressure routing. A novel feature of our scheme is that a utility function can be assigned to a collection of flows, generalizing the traditional approach of optimizing per-flow utilities. Our control policies use finite-buffer queues and are independent of arrival statistics. Their near-optimal performance is proved and further supported by simulation results.
\end{abstract}

\section{Introduction}

The idea behind flow control in data networks is to regulate the source rates in order to prevent network overload, and provide fair allocation of resources. In recent years, extensive research has been devoted to the problem of network utility maximization, with the objective of optimizing network resource allocation through a combination of source-based flow control, routing, and scheduling.  The common theme is to assign a utility, as a function of the source rate, to a flow specified by a source-destination pair, and formulate an optimization problem that maximizes the sum of utilities (e.g., see~\cite{CLC07,KMT98,Kel97,LaL99,Sto05,NML08,EaS06a,EaS07,LaS04conf}). The optimal network control policy is revealed as the algorithm that solves the utility maximization problem. 


Source-based flow control implicitly requires all sources to react properly to congestion signals such as packet loss or delay. Congestion-insensitive traffic, however, is pervasive in modern networks. UDP-based applications such as video streaming or VoIP are increasingly popular and do not respond to congestion. Moreover, greedy or malicious users can inject excessive data into the network to either boost self-performance or bring down high-profile websites. In these circumstances, the network can be temporarily overloaded and congestion-aware applications, e.g., TCP-based traffic, may be adversely affected or even starved. In this context, source-based flow control may not be adequate.

There are other scenarios in which optimizing the source rates of data flows on a per-flow basis is ineffective. The Internet nowadays is vulnerable to distributed denial of service (DDoS) attacks~\cite{Cha02,SGT11}, in which an attacker creates a large number of flows with different source addresses to overwhelm prominent websites. Similarly, in a multi-player online game, thousands of users require continuous access to the game servers. There are also occasions in which a large number of users seek access to a website that broadcasts live events~\cite{CNET11news}; this is the so-called \emph{flash crowd} phenomenon~\cite{MBF02}.  In these situations, source-based flow control is ineffective because each individual flow uses little bandwidth, but their aggregate traffic can lead to severe link congestion near the service provider, and may starve other users in the network. A flow control scheme that can  optimize a utility assigned to a collection of flows as opposed to optimizing the sum of per-flow utilities can be used to cope with such scenarios. Moreover, in order to cope with uncooperative flows, it is necessary to relocate the flow control functionality from untrusted network hosts to more secure ones, such as the web servers that provide service at the receiver end. 

\begin{figure}[htbp]
\centering
\includegraphics[width=3.5in]{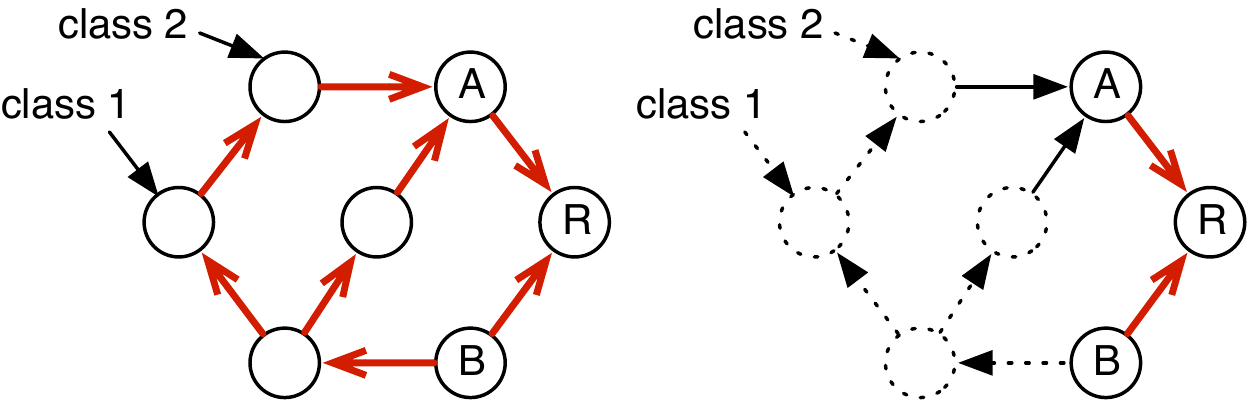}
\caption{Our receiver-based policy can be implemented in the whole network on the left, or implemented only at nodes $A$, $B$, and $R$ on the right, where $R$ is the only receiver. The rest of the network on the right may be controlled by another network operator or follow a different network control scheme.}
\label{fig:105}
\end{figure}

In this paper, we develop such \emph{receiver-based flow control} policies using tools from stochastic network optimization~\cite{GNT06,Nee10book}. Our main contributions are three-fold. First, we formulate a utility maximization problem that can assign utilities to an aggregate of flows, of which the usual per-flow-based utility maximization is a special case. Second, given an arbitrary arrival rate matrix (possibly outside the network's stability region), we characterize the corresponding achievable throughput region in terms of queue overflow rates. Third, using a novel decomposition of the utility functions, we design a network control policy consisting of: (i) a set of flow controllers at the receivers; (ii) packet dropping mechanism at internal nodes; and (iii) back-pressure routing at intermediate nodes. The receiver-based flow controllers adjust throughput by modifying the differential backlogs between the receivers and their neighboring nodes---a small (or negative) differential backlog is regarded as a \emph{push-back} mechanism to slow down data delivery to the receiver. To deal with data that cannot be delivered due to network overload, we design a threshold-based packet dropping mechanism that discards data whenever queues grow beyond certain thresholds. Surprisingly, we show that our threshold-based packet dropping scheme, without the use of flow control, is sufficient to maximize the weighted sum throughput. Moreover, the combined flow control and packet dropping mechanism has the following properties: (i) It works with finite-size buffers. (ii) It is nearly utility-optimal (throughput-optimal as a special case) and the performance gap from the optimal utility goes to zero as buffer sizes increase. (iii) It does not require the knowledge of arrival rates and therefore is robust to time-varying arrival rates that can go far beyond the network's stability region. In addition, our control policy can be implemented only in parts of the network that include the receivers, treating the rest of the network as exogenous data sources (see Fig.~\ref{fig:105} for an example), and thus might be an attractive flow control solution for web servers.



There has been a significant amount of research in the general area of stochastic network control.  Utility-optimal policies that combine source-end flow control with back-pressure routing have been studied in~\cite{NML08,EaS06a,EaS07,LaS04conf} (and references therein). These policies optimize per-flow utilities and require infinite-capacity buffers. However, they are not robust in the face of uncooperative users who may not adhere to the flow control scheme. A closely related problem to that studied in this paper is that of characterizing the queue overflow rates in lossless networks in overload.  In a single-commodity network, a back-pressure policy is shown to achieve the most balanced queue overflow rates~\cite{GaT06}, and controlling queue growth rates using the max-weight policy is discussed in~\cite{CAB11draft}. The queue growth rates in networks under max-weight and $\alpha$-fairness policies are analyzed in~\cite{SaW11, EBZ08conf}. We finally note that the importance of controlling an aggregate of data flows has been  addressed in~\cite{MBF02}, and rate-limiting mechanisms in front of a web server to achieve some notion of max-min fairness against DDoS attacks have been proposed in~\cite{YLL05,TCL07, CaS05}.

An outline of the paper is as follows. The network model is given in Section~\ref{sec:601}. We formulate the utility maximization problem and characterize the achievable throughput region in terms of queue overflow rates in Section~\ref{sec:602}. Section~\ref{sec:801} introduces a threshold-based packet dropping policy that maximizes the weighted sum throughput without the use of flow control. Section~\ref{sec:501} presents a receiver-based flow control and packet dropping policy that solves the general utility maximization problem. Simulation results that demonstrate the near-optimal performance of our policies are given in Sections~\ref{sec:801} and~\ref{sec:501}.

\section{Network Model} \label{sec:601}
We consider a network with nodes $\cN = \{1, 2, \ldots, N\}$ and directed links $\cL = \{(n,m) \mid n,m\in \cN\}$.  Assume time is slotted. In every slot, packets randomly arrive at the network for service and are categorized into a collection $\cC$ of classes. The definition of a data class is quite flexible except that we assume packets in a class $c \in \cC$ have a shared destination~$d_{c}$. For example, each class can simply be a flow specified by a source-destination pair. Alternatively, computing-on-demand services in the cloud such as Amazon (Elastic Compute Cloud; EC2) or Google (App Engine) can assign business users to one class and residential users to another. Media streaming applications may categorize users into classes according to different levels of subscription to the service provider. While classification of users/flows in various contexts is a subject of significant importance, in this paper we assume for simplicity that the class to which a packet belongs can be ascertained from information contained in the packet (e.g., source/destination address, tag, priority field, etc.).
Let $A_{n}^{(c)}(t)\in \{0, 1, \ldots, \amax\}$ be the number of exogenous class~$c$ packets arriving at node $n$ in slot $t$, where $\amax$ is a finite constant; let $A_{d_{c}}^{(c)}(t)=0$ for all $t$. We assume $A_{n}^{(c)}(t)$ are independent across classes $c$ and nodes $n\neq d_{c}$, and are i.i.d. over slots with mean $\mathbb{E}\big[A_{n}^{(c)}(t)\big] = \lambdanc$. 

In the network, packets are relayed toward the destinations via dynamic routing and link rate allocation decisions. Each link $(n,m)\in\cL$ is used to transmits data from node $n$ to node $m$ and has a fixed capacity $\mu_{\text{max}}^{nm}$ (in units of packets/slot).\footnote{We focus on wireline networks in this paper for ease of exposition. Our results and analysis can be easily generalized to wireless networks or switched networks in which link rate allocations are subject to interference constraints.} Under a given control policy, let $\mu_{nm}^{(c)}(t)$ be the service rate allocated to class~$c$ data over link $(n,m)$ in slot $t$. The service rates must satisfy the link capacity constraints
\[
\sum_{c\in \cC} \mu_{nm}^{(c)}(t) \leq \mu_{\text{max}}^{nm}, \quad \text{for all $t$ and all links $(n,m)$.}
\]

At a node $n\neq d_{c}$, class $c$ packets that arrive but not yet transmitted are stored in a queue; we let $\Qnc(t)$ be the backlog of class~$c$ packets at node $n$ at time $t$. We assume initially $\Qnc(0) = 0$ for all $n$ and $c$. Destinations do not buffer packets and we have $Q_{d_{c}}^{(c)}(t) = 0$ for all $c$ and $t$. For now, we assume every queue $\Qnc(t)$ has an infinite-capacity buffer; we show later that our control policy needs only finite buffers. To resolve potential network congestion due to traffic overload, a queue $\Qnc(t)$, after transmitting data to neighboring nodes in a slot, discards $d_{n}^{(c)}(t)$ packets from the remaining backlog at the end of the slot. The drop rate $d_{n}^{(c)}(t)$ is a function of the control policy to be described later and takes values in $[0, \dmax]$ for some finite $\dmax$. The queue $\Qnc(t)$ evolves over slots according to
\begin{equation} \label{eq:201}
\begin{split}
\Qnc(t+1) &\leq \bigg[ \Big( \Qnc(t) - \sum_{b} \mu_{nb}^{(c)}(t) \Big)^{+} - \dnc(t) \bigg]^{+} \\
&\quad+ \Anc(t) + \sum_{a} \mu_{an}^{(c)}(t), \quad \forall c,\ \forall n\neq d_{c},
\end{split}
\end{equation}
where $(\cdot)^{+} \triangleq \max(\cdot, 0)$.  This inequality  is due to the fact that endogenous arrivals may be less than the allocated rate $\sum_{a} \mu_{an}^{(c)}(t)$ when neighboring nodes do not have sufficient packets to send. 

For convenience,  we define the maximum transmission rate into and out of a node by
\begin{equation} \label{eq:1204}
\mumaxin \triangleq \max_{n\in \cN} \sum_{a:(a,n)\in\cL} \mumax^{an}, \quad \mumaxout \triangleq \max_{n\in \cN} \sum_{b: (n,b)\in\cL} \mumax^{nb}.
\end{equation}
Throughout the paper, we use the following assumption.

\begin{assu} \label{assu:601}
We assume $\dmax \geq \amax + \mumaxin$.
\end{assu}
From~\eqref{eq:201}, the sum $\amax+\mumaxin$ is the largest amount of data that can arrive at a node in a slot; therefore it is an upper bound on the maximum queue overflow rate at any node. Assumption~\ref{assu:601} ensures that the maximum packet dropping rate $\dmax$ is no less than the maximum queue overflow rate, so that we can always prevent the queues from blowing up.

\section{Problem Formulation} \label{sec:602}

We assign to each class $c\in\cC$ a utility function $g_{c}(\cdot)$. Given an (unknown) arrival rate matrix $\blambda = (\lambdanc)$, let $\Lambda_{\blambda}$ be the set of all achievable throughput vectors $(r_{c})$, where $r_{c}$ is the aggregate throughput of class $c$ received by the destination $d_{c}$. Note that $\Lambda_{\blambda}$ is a function of~$\blambda$. We seek to design a control policy that solves the global utility maximization problem
\begin{align}
\text{maximize} &\quad \sum_{c\in\cC} g_{c}(r_{c}) \label{eq:207} \\
\text{subject to} &\quad (r_{c}) \in \Lambda_{\blambda}, \label{eq:208}
\end{align}
where the region $\Lambda_{\blambda}$ is presented later in Lemma~\ref{lem:201}. We assume all $g_{c}(\cdot)$ functions are concave, increasing, and continuously differentiable. For ease of exposition, we also assume the functions $g_{c}(\cdot)$ have bounded derivatives such that $\abs{g_{c}'(x)} \leq m_{c}$ for all $x\geq 0$, where $m_{c}$ are finite constants.\footnote{Utility functions $g_{c}(x)$ that have unbounded derivatives as $x\to 0$ can be approximated by those with bounded derivatives. For example, we may approximate the proportional fairness function $\log(x)$ by $\log(x+\xi)$ for some small $\xi >0$.} 

\begin{figure}[htbp]
\centering
\includegraphics[width=2in]{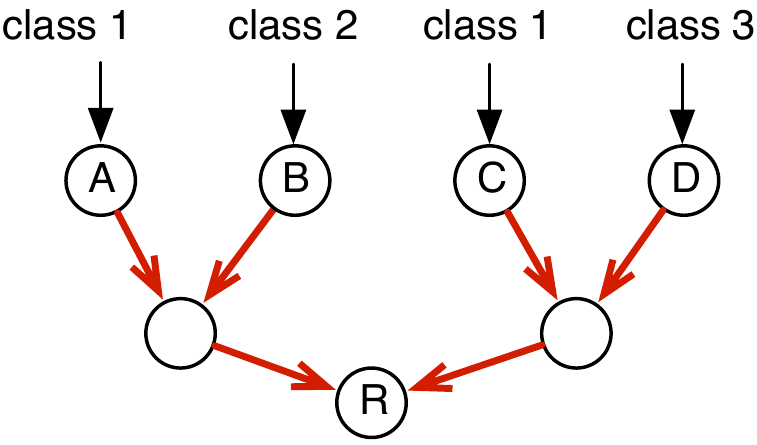}
\caption{A tree network with three classes of traffic.}
\label{fig:102}
\end{figure}

As an example, consider the tree network in Fig.~\ref{fig:102} that serves three classes of traffic destined for node $R$.
Class~$1$ data originates from two different sources $A$ and $C$, and may represent the collection of users located in different parts of the network sending or requesting information from node $R$.  If class $1$ traffic is congestion-insensitive and overloads the network, without proper flow control class $2$ and $3$ will be starved due to the presence of class $1$. A utility maximization problem here is to solve
\begin{align}
\text{maximize} &\quad g_{1}(r_{1A}+r_{1C}) + g_{2}(r_{2B}) + g_{3}(r_{3D}) \label{eq:1201} \\
\text{subject to} &\quad \text{$(r_{1A}, r_{1C}, r_{2B}, r_{3D})$ feasible,} \label{eq:1202}
\end{align}
where $r_{1A}$ denotes the throughput of class $1$ data originating from $A$;  $r_{1C}$, $r_{2B}$, and $r_{3D}$ are  defined similarly. Note that this utility maximization~\eqref{eq:1201}-\eqref{eq:1202} is very different from, and generalizes, the traditional per-flow-based utility maximization.

\subsection{Achievable Throughput Region}

The next lemma characterizes the set $\Lambda_{\blambda}$ of all achievable throughput vectors in~\eqref{eq:208}.

\begin{lem} \label{lem:201}
Under i.i.d. arrival processes with an arrival rate matrix $ \blambda = (\lambdanc)$, let $\Lambda_{\blambda}$ be the closure of the set of all achievable throughput vectors $(r_{c})$. Then $(r_{c}) \in \Lambda_{\blambda}$ if and only if there exist  flow variables $\{f_{ab}^{(c)} \geq 0\mid c\in\cC, (a,b)\in \cL\}$ and queue overflow variables $\{\qnc\geq 0 \mid c\in \cC,\ n\neq d_{c}\}$ such that
\begin{gather}
\lambdanc + \sum_{a} f_{an}^{(c)} = q_{n}^{(c)} + \sum_{b} f_{nb}^{(c)} \quad \forall\, c\in \cC,\ n\neq d_{c}, \label{eq:202} \\
\sum_{c} f_{ab}^{(c)} \leq \mumax^{ab} \quad \forall\, (a,b)\in \cL, \label{eq:203}\\
r_{c} \leq \sum_{a} f_{ad_{c}}^{(c)} = \sum_{n} \lambdanc - \sum_{n} \qnc \quad\forall\, c\in \cC. \label{eq:204}
\end{gather}
In other words,
\[
\Lambda_{\blambda} = \Set{ (r_{c}) | \text{\eqref{eq:202}-\eqref{eq:204} hold and } f_{ab}^{(c)} \geq 0, \qnc \geq 0}.
\]
\end{lem}

\begin{IEEEproof}[Proof of Lemma~\ref{lem:201}]
See Appendix~\ref{appendix:201}.
\end{IEEEproof}

In Lemma~\ref{lem:201}, equation~\eqref{eq:202} is the flow conservation constraint stating that the total flow rate of a class into a node is equal to the flow rate out of the node plus the queue overflow rate. Equation~\eqref{eq:203} is the link capacity constraint. The equality in~\eqref{eq:204} shows that the throughput of a class is equal to the sum of exogenous arrival rates less the queue overflow rates. Lemma~\ref{lem:201} is closely related to the network capacity region $\Lambda$ defined in terms of admissible arrival rates (see Definition~\ref{defn:101}); their relationship is shown in the next corollary.

\begin{defn}[Theorem~$1$,~\cite{NMR05}] \label{defn:101}
The capacity region $\Lambda$ of the network is the set of all arrival rate matrices $(\lambda_{n}^{(c)})$ for which there exists nonnegative flow variables $f_{ab}^{(c)}$ such that
\begin{gather}
\lambda_{n}^{(c)} + \sum_{a} f_{an}^{(c)} \leq \sum_{b} f_{nb}^{(c)}, \quad \forall\,c\in \cC,\ n\neq d_{c}, \label{eq:112} \\
\sum_{c} f_{ab}^{(c)} \leq \mumax^{ab}, \quad \forall\, (a,b)\in  \cL. \label{eq:113}
\end{gather}
\end{defn}

\begin{cor} \label{lem:501}
An arrival rate matrix $(\lambdanc)$ lies in the network capacity region $\Lambda$ if there exist flow variables $f_{ab}^{(c)} \geq 0$ such that flow conservation constraints~\eqref{eq:202} and link capacity constraints~\eqref{eq:203} hold with $\qnc = 0$ for all $n$ and $c$.
\end{cor}

Corollary~\ref{lem:501} shows that $(\lambdanc)$ is achievable if and only if there exists a control policy yielding zero queue overflow rates. In this case the throughput of class $c$ is $r_{c} = \sum_{n} \lambdanc$.

We remark that the solution to the utility maximization~\eqref{eq:207}-\eqref{eq:208} activates the linear constraints~\eqref{eq:204}; thus the problem~\eqref{eq:207}-\eqref{eq:208} is equivalent to
\begin{align}
\text{maximize} &\quad \sum_{c\in\cC} g_{c}(r_{c}) \label{eq:212} \\
\text{subject to} &\quad r_{c} = \sum_{n} \lambdanc - \sum_{n} \qnc, \quad \forall \, c\in \cC \label{eq:213} \\
&\quad \text{\eqref{eq:202} and~\eqref{eq:203} hold} \\
&\quad f_{ab}^{(c)} \geq 0,\ q_{n}^{(c)}\geq 0,\ \forall (a,b)\in\cL,\ \forall n,c. \label{eq:605}
\end{align}
Let $(r_{c}^{*})$ be the optimal throughput vector that solves~\eqref{eq:212}-\eqref{eq:605}. If the arrival rate matrix $(\lambdanc)$ is in the network capacity region $\Lambda$, the optimal throughput is $r_{c}^{*} = \sum_{n} \lambdanc$ from Corollary~\ref{lem:501}. Otherwise, we have $r_{c}^{*} = \sum_{n}\lambdanc - \sum_{n} q_{n}^{(c)*}$, where $q_{n}^{(c)*}$ is the optimal queue overflow rate.

\subsection{Features of Network Control} \label{sec:603}

Our control policy that solves~\eqref{eq:212}-\eqref{eq:605} has two main features. First, we have a packet dropping mechanism discarding data from the network when queues build up. An observation here is that, in order to optimize throughput and keep the network stable, we should drive the packet dropping rate to be equal to the optimal queue overflow rate. Second, we need a flow controller driving the throughput vector toward the utility-optimal point. To convert the control objective~\eqref{eq:212} into these two control features, we define, for each class $c\in\cC$, a utility function $h_{c}(\cdot)$ related to $g_{c}(\cdot)$ as
\begin{equation} \label{eq:606}
h_{c}(r_{c}) \triangleq g_{c}(r_{c}) - \theta_{c} \, r_{c},
\end{equation}
where $\theta_{c} \geq 0$ are control parameters to be decided later. Using~\eqref{eq:213}, we have
\[
g_{c}(r_{c}) = h_{c}(r_{c}) + \theta_{c} \Big[ \sum_{n} \lambdanc - \sum_{n} \qnc \Big].
\]
Since $\lambdanc$ are unknown constants, maximizing $\sum_{c\in\cC} g_{c}(r_{c})$ is the same as maximizing
\begin{equation} \label{eq:607}
\sum_{c\in\cC} \Big[ h_{c}(r_{c}) -\theta_{c} \sum_{n} \qnc \Big].
\end{equation}
This equivalent objective~\eqref{eq:607} can be optimized by jointly maximizing the new utility $\sum_{c\in\cC} h_{c}(r_{c})$ at the receivers and minimizing the weighted queue overflow rates (i.e., the weighted packet dropping rates) $\sum_{c\in\cC} \theta_{c}\, q_{n}^{(c)}$ at each node $n$.

Optimizing the throughput at the receivers amounts to controlling the amount of data \emph{actually delivered}.  This is difficult because the data available to the receivers at their upstream nodes is highly correlated with control decisions taken in the rest of the network. Optimizing the packet dropping rates depends on the data available at each network node, which has similar difficulties. To get around these difficulties, we introduce auxiliary control variables $\varphinc\geq 0$ and $\nu_{c}\geq 0$ and consider the optimization problem
\begin{align}
\text{maximize} &\quad \sum_{c} \Big[ h_{c}(\nu_{c}) - \theta_{c} \sum_{n} \varphinc \Big] \label{eq:901} \\
\text{subject to} &\quad r_{c} = \nu_{c}, \quad \forall\, c\in \cC \label{eq:210}, \\
&\quad \qnc \leq \varphinc, \quad \forall\, c\in \cC,\ n\neq d_{c}. \label{eq:211} \\
&\quad \text{\eqref{eq:213}-\eqref{eq:605} hold.} \label{eq:902}
\end{align}
This is an equivalent problem to~\eqref{eq:212}-\eqref{eq:605}. The constraints~\eqref{eq:210} and~\eqref{eq:211} can be enforced by stabilizing virtual queues that will appear in our control policy. The new control variables $\nu_{c}$ and $\varphinc$  to be optimized can now be chosen freely unconstrained by past control actions in the network. Introducing auxiliary variables and setting up virtual queues are at the heart of using Lyapunov drift theory to solve network optimization problems.

\section{Maximizing the Weighted Sum Throughput} \label{sec:801}

For ease of exposition, we first consider the special case of maximizing the weighted sum throughput in the network. For each class $c\in\cC$, we let $g_{c}(r_{c}) = a_{c}\,r_{c}$ for some $a_{c} > 0$. We present a threshold-based packet dropping policy that, together with back-pressure routing, solves this problem. Surprisingly, flow control is not needed here. This is because maximizing the weighted sum throughput is equivalent to minimizing the weighted packet dropping rate. Indeed, choosing $\theta_{c} =  a_{c}$  in~\eqref{eq:606}, we have $h_{c}=0$ for all classes $c$, under which maximizing the equivalent objective~\eqref{eq:607} is the same as minimizing $\sum_{n,c} \theta_{c}\, \qnc$. In the next section, we will combine the threshold-based packet dropping policy with receiver-based flow control to solve the general utility maximization problem. 

\subsection{Control Policy}

To optimize packet dropping rates, we set up a \emph{drop queue} $\Dnc(t)$ associated with each queue $\Qnc(t)$. The packets that are dropped from $\Qnc(t)$ in a slot, denoted by $\tdnc(t)$, are first stored in $\Dnc(t)$ for eventual deletion. From~\eqref{eq:201}, we have
\begin{equation} \label{eq:1203}
\tdnc(t) = \min \bigg[ \Big(\Qnc(t) - \sum_{b} \mu_{nb}^{(c)}(t) \Big)^{+},\ \dnc(t) \bigg].
\end{equation}
Note that the quantity $\tdnc(t)$ is the \emph{actual} packets dropped from $\Qnc(t)$, which is strictly less than the allocated drop rate $\dnc(t)$ if queue $\Qnc(t)$ does not have sufficient data. Packets are permanently deleted from $\Dnc(t)$ at the rate of $\varphinc(t) \in [0, \dmax]$ in slot $t$. The queue $\Dnc(t)$ evolves according to
\begin{equation} \label{eq:229}
\Dnc(t+1) = \big[ \Dnc(t) - \varphinc(t) \big]^{+}+ \tdnc(t).
\end{equation}
Assume initially $\Dnc(0) = V\theta_{c} = V a_{c}$ for all $n$ and $c$, where $V>0$ is a control parameter.\footnote{It suffices to assume $\Dnc(0)$ to be finite. Our choice of $\Dnc(0) = V\theta_{c}$ avoids unnecessary packet dropping in the initial phase of the system.} If queue $\Dnc(t)$ is stabilized, then minimizing the service rate of $\Dnc(t)$ effectively minimizes the time average of dropped packets  at $\Qnc(t)$. We propose the following policy.

\noindent \rule[0.05in]{3.5in}{0.01in}

Overload Resilient Algorithm ($\ORA$)

\noindent \rule[0.05in]{3.5in}{0.01in}

\emph{Parameter Selection:} Choose $\theta_{c} = a_{c}$ for all classes $c\in\cC$, where $g_{c}(x) = a_{c}\, x$. Choose a parameter $V>0$.

\emph{Backpressure Routing:} Over each link $l = (n,m)\in \cL$, let $\cC_{l}$ be the subset of classes that have access to link $l$. Compute the differential backlog $W_{l}^{(c)}(t) = Q_{n}^{(c)}(t) - Q_{m}^{(c)}(t)$ for each class $c\in \cC_{l}$, where  $Q_{d_{c}}^{(c)}(t)=0$ at the receiver $d_{c}$. Define
\begin{gather*}
W_{l}^{(c)*}(t) = \max_{c\in \cC_{l}} W_{l}^{(c)}(t), \\
c_{l}^{*}(t) = \argmax_{c\in \cC_{l}} W_{l}^{(c)}(t).
\end{gather*}
We allocate the service rates
\[
\mu_{nm}^{(c_{l}^{*}(t))}(t) = \begin{cases} \mumax^{nm} & \text{if $W_{l}^{(c)*}(t) > 0,$} \\ 0 & \text{if $W_{l}^{(c)*}(t) \leq 0$.}\end{cases}
\]
Let $\mu_{nm}^{(c)}(t)=0$ for all classes $c = \cC_{l} \setminus \{c_{l}^{*}(t)\}$.

\emph{Packet Dropping:} At queue $\Qnc(t)$, allocate the packet dropping rate $\dnc(t)$ (see~\eqref{eq:201}) according to
\[
\dnc(t) = \begin{cases}
\dmax & \text{if $\Qnc(t) > \Dnc(t)$,} \\ 0 & \text{if $\Qnc(t) \leq \Dnc(t)$,}
\end{cases}
\]
where $\dmax>0$ is a constant chosen to satisfy Assumption~\ref{assu:601}. At the drop queue $\Dnc(t)$, allocate its service rate $\varphinc(t)$ according to
\[
\varphinc(t) = \begin{cases}
\dmax & \text{if $\Dnc(t) > V\theta_{c}$,} \\ 0 & \text{if $\Dnc(t) \leq V\theta_{c}$.}
\end{cases}
\]

\emph{Queue Update:} Update queues $\Qnc(t)$ according to~\eqref{eq:201} and update queues $\Dnc(t)$ according to~\eqref{eq:1203}-\eqref{eq:229} in every slot.

\noindent \rule[0.05in]{3.5in}{0.01in}

The packet dropping subroutine in this policy is threshold-based. The $\ORA$ policy uses local queueing information and does not require the knowledge of exogenous arrival rates. It is notable that network overload is autonomously resolved by each node making local decisions of routing and packet dropping.

\subsection{Performance of the $\ORA$ Policy}

\begin{lem}[Deterministic Bound for Queues] \label{lem:402}
For each class $c\in\cC$, define the constants
\begin{equation} \label{eq:503}
\Dmax^{(c)} \triangleq V \theta_{c} + \dmax, \quad \Qmax^{(c)} \triangleq V\theta_{c} + 2\dmax.
\end{equation}
In the $\ORA$ policy, queues $\Qnc(t)$ and $\Dnc(t)$ are deterministically bounded by
\[
\Qnc(t)\leq \Qmax^{(c)},\quad \Dnc(t)\leq \Dmax^{(c)}, \quad\text{for all $t$, $c$, and $n\neq d_{c}$.}
\]
In addition, we have $\Dnc(t) \geq V\theta_{c}  - \dmax$ for all $n$, $c$, and $t$.
\end{lem}

\begin{IEEEproof}
See Appendix~\ref{appendix:601}.
\end{IEEEproof}

In Lemma~\ref{lem:402}, the value of $\Qmax^{(c)}$ is the finite buffer size sufficient at queue $\Qnc(t)$. The parameter $V$ controls when queue $\Qnc(t)$ starts dropping packets. Indeed, due to $\Dnc(t) \geq V\theta_{c} - \dmax$, the $\ORA$ policy discards packets from $\Qnc(t)$ only if $\Qnc(t) \geq V\theta_{c}-\dmax$. The quantity $V\theta_{c}-\dmax$ is a controllable threshold beyond which we say queue $\Qnc(t)$ is overloaded and should start dropping packets. As we see next, the performance of the $\ORA$ policy approaches optimality as the buffer sizes increase.

\begin{thm} \label{thm:501}
Define the limiting throughput of class $c$ as
\begin{equation} \label{eq:504}
\overline{r_{c}} \triangleq \lim_{t\to\infty} \frac{1}{t} \sum_{\tau=0}^{t-1} \mathbb{E}\bigg[ \sum_{a: (a, d_{c})\in\cL} \widetilde{\mu}_{ad_{c}}^{(c)}(\tau)\bigg],
\end{equation}
where $\widetilde{\mu}_{ad_{c}}^{(c)}(\tau)$ denotes the class $c$ packets received by node $d_{c}$ over link $(a, d_{c})$. The $\ORA$ policy yields the limiting weighted sum throughput satisfying
\begin{equation} \label{eq:505}
\sum_{c\in\cC} a_{c}\, \overline{r_{c}} \geq \sum_{c\in\cC} a_{c} \, r_{c}^{*} - \frac{B}{V},
\end{equation}
where $(r_{c}^{*})$ is the optimal throughput vector that solves~\eqref{eq:212}-\eqref{eq:605} under the linear objective function $\sum_{c\in\cC} a_{c}\,r_{c}$, $V>0$ is a control parameter, and $B$ is a finite constant defined as
\[
B \triangleq \abs{\cN}\abs{\cC} \left[(\mumaxout + \dmax)^{2} + (\amax + \mumaxin)^{2} + 2\dmax^{2}\right],
\]
where $\abs{\mathcal{A}}$ denotes the cardinality of a set $\mathcal{A}$.
\end{thm}
We omit the proof of Theorem~\ref{thm:501} because it is similar to that of Theorem~\ref{thm:502} presented later in the general case of utility maximization.  From~\eqref{eq:505}, the $\ORA$  policy yields near-optimal performance by choosing the parameter $V$ sufficiently large. Correspondingly, a large $V$ implies a large buffer size of $\Qmax^{(c)} = V\theta_{c} + 2\dmax$.

As shown in Corollary~\ref{lem:501}, if the arrival rate matrix $(\lambdanc)$ lies in the network capacity region $\Lambda$, then the optimal throughput for class $c$ is $r_{c}^{*} = \sum_{n} \lambdanc$ and~\eqref{eq:505} reduces to
\[
\sum_{c\in\cC} a_{c}\, \overline{r_{c}} \geq \sum_{c\in\cC} a_{c} \Big(\sum_{n} \lambdanc\Big) - \frac{B}{V}.
\]
That we can choose $V$ arbitrarily large leads to the next corollary.
\begin{cor} \label{cor:601}
The $\ORA$ policy is (close to) throughput optimal.
\end{cor}

\subsection{Simulation of the $\ORA$ Policy}

\begin{figure}[htbp]
\centering
\includegraphics[width=2in]{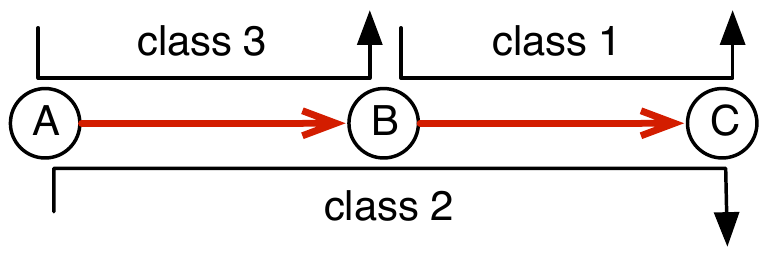}
\caption{A 3-node network with three classes of traffic.}
\label{fig:101}
\end{figure}

We conduct simulations for the $\ORA$ policy in the network shown in Fig.~\ref{fig:101}. 
The directed links $(A, B)$ and $(B, C)$ have the capacity of $1$ packet/slot. There are three classes of traffic to be served; for example, class 1 data arrives at node $B$ and is destined for node $C$. Classes $1$ and $2$ compete for service over $(B, C)$; classes $2$ and $3$ compete for service over $(A, B)$. Each simulation below is run over $10^{6}$ slots.

\subsubsection{Fixed arrival rates} In each class, we assume a Bernoulli arrival process whereby $20$ packets arrive to the network in a slot with probability $0.1$, and no packets arrive otherwise.  The arrival rate of each class is 2 packets/slot, which clearly overloads the network. 

Let $r_{c}$ be the throughput of class $c$. Consider the objective of maximizing the weighted sum throughput $3\, r_{1} +  2\, r_{2} + r_{3}$; the weights are rewards obtained by serving a packet in a class. The optimal solution is: (i) Always serve class $1$ at node $B$ because it yields better rewards than serving class $2$. (ii) Always serve class $3$  at node $A$---although class $2$ has better rewards than class $3$, it does not make sense to serve class $2$ at $A$ only to be dropped later at $B$. The optimal throughput vector is therefore $(1, 0, 1)$.  Consider another objective of maximizing $3\, r_{1} +  5\, r_{2} + r_{3}$. Here, class $2$ has a reward that is better than the sum of rewards of the other two classes. Thus both nodes $A$ and $B$ should always serve class $2$; the optimal throughput vector  is $(0, 1, 0)$. Table~\ref{table:601} shows the near-optimal performance of the $\ORA$ policy in both cases as $V$ increases.
\begin{table}[tp]
\caption{The throughput performance of the $\ORA$ policy under fixed arrival rates.}
\label{table:601}
\centering
\small
\subtable[Maximizing $3\,r_{1}+2\,r_{2}+r_{3}$]{
\begin{tabular}{rccc}
\toprule
$V$ & $r_{1}$ & $r_{2}$ & $r_{3}$ \\ \otoprule
$10$ & $.787$  & $.168$ & $.099$  \\
$20$ & $.867$  & $.133$ & $.410$ \\
$50$ & $.992$  & $.008$ & $.967$ \\
$100$ & $.999$ & $0$ & $.999$ \\  \otoprule
opt & $1$ & $0$ & $1$ \\
\bottomrule
\end{tabular}
}
\subtable[Maximizing $3\,r_{1}+5\,r_{2}+r_{3}$]{
\begin{tabular}{rccc}
\toprule
$V$ & $r_{1}$ & $r_{2}$ & $r_{3}$ \\ \otoprule
$10$ & $.185$ & $.815$ & $.083$ \\
$20$ & $.107$ & $.893$ & $.095$ \\
$50$ & $.031$ & $.969$ & $.031$ \\
$100$ & $.002$ & $.998$ & $.001$ \\ \otoprule
opt & $0$ & $1$ & $0$ \\
\bottomrule
\end{tabular}
}
\end{table}



\subsubsection{Time-varying arrival rates}
We show that the $\ORA$ policy is robust to time-varying arrival rates. Suppose class $1$ and $3$ have a fixed arrival rate of $0.8$ packets/slot. The arrival rate of class $2$ is $2$ packets/slot in the  interval $\mathcal{T} = [3\times10^{5}, 6\times 10^{5})$ and is $0.1$ packets/slot elsewhere.  We consider the objective of maximizing $3\, r_{1} +  5\, r_{2} + r_{3}$. The network is temporarily overloaded in the interval $\mathcal{T}$; the optimal time-average throughput in $\mathcal{T}$ is $(0, 1, 0)$ as explained in the above case. The network is underloaded in the interval $[0, 10^{6}) \setminus \mathcal{T}$, in which the optimal throughput is $(0.8, 0.1, 0.8)$.

\begin{figure}[tp]
\centering 
\subfigure[Queue $Q_{B}^{(1)}(t)$]{\includegraphics[width=1.7in]{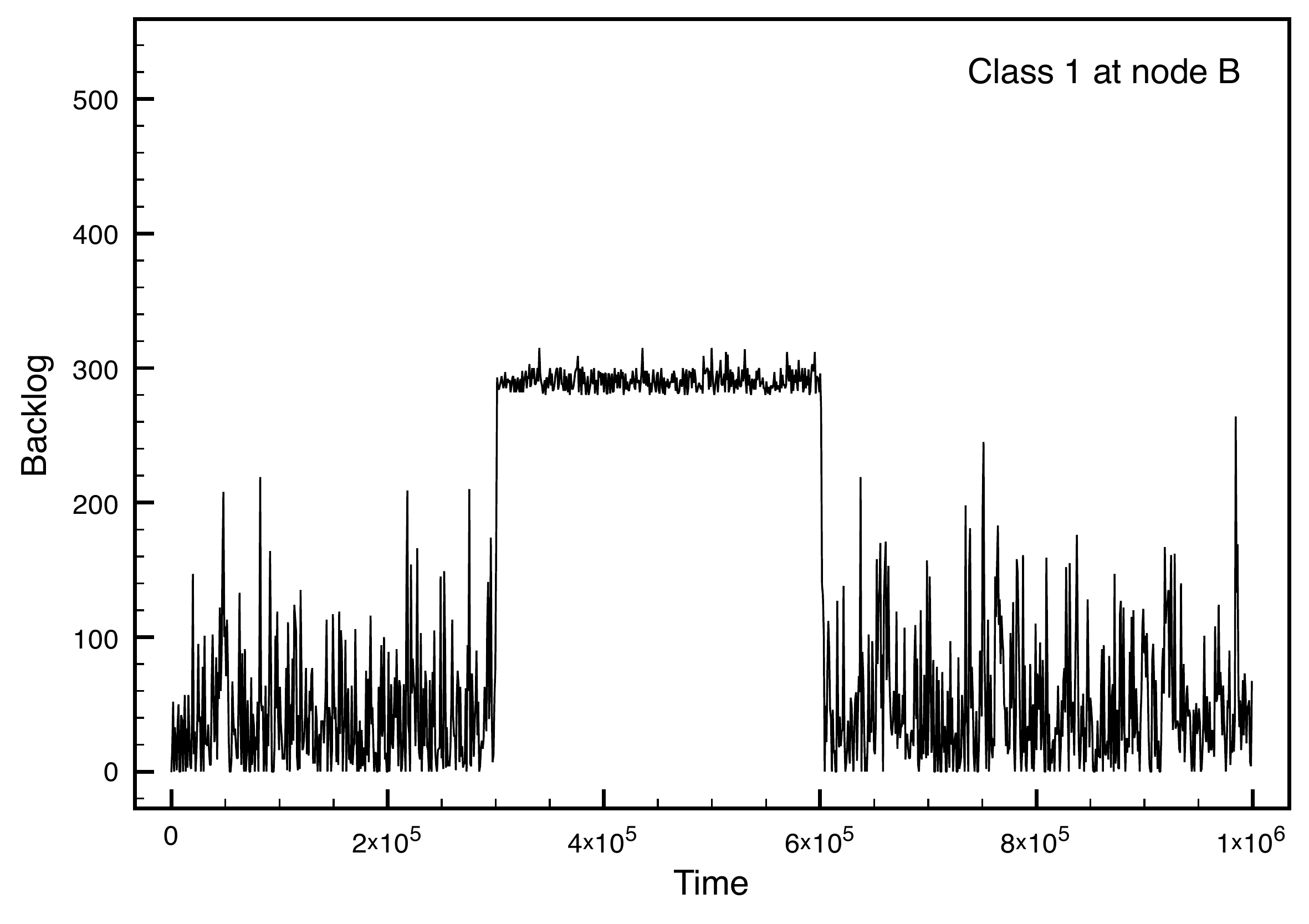}}
\subfigure[Queue $Q_{B}^{(2)}(t)$]{\includegraphics[width=1.7in]{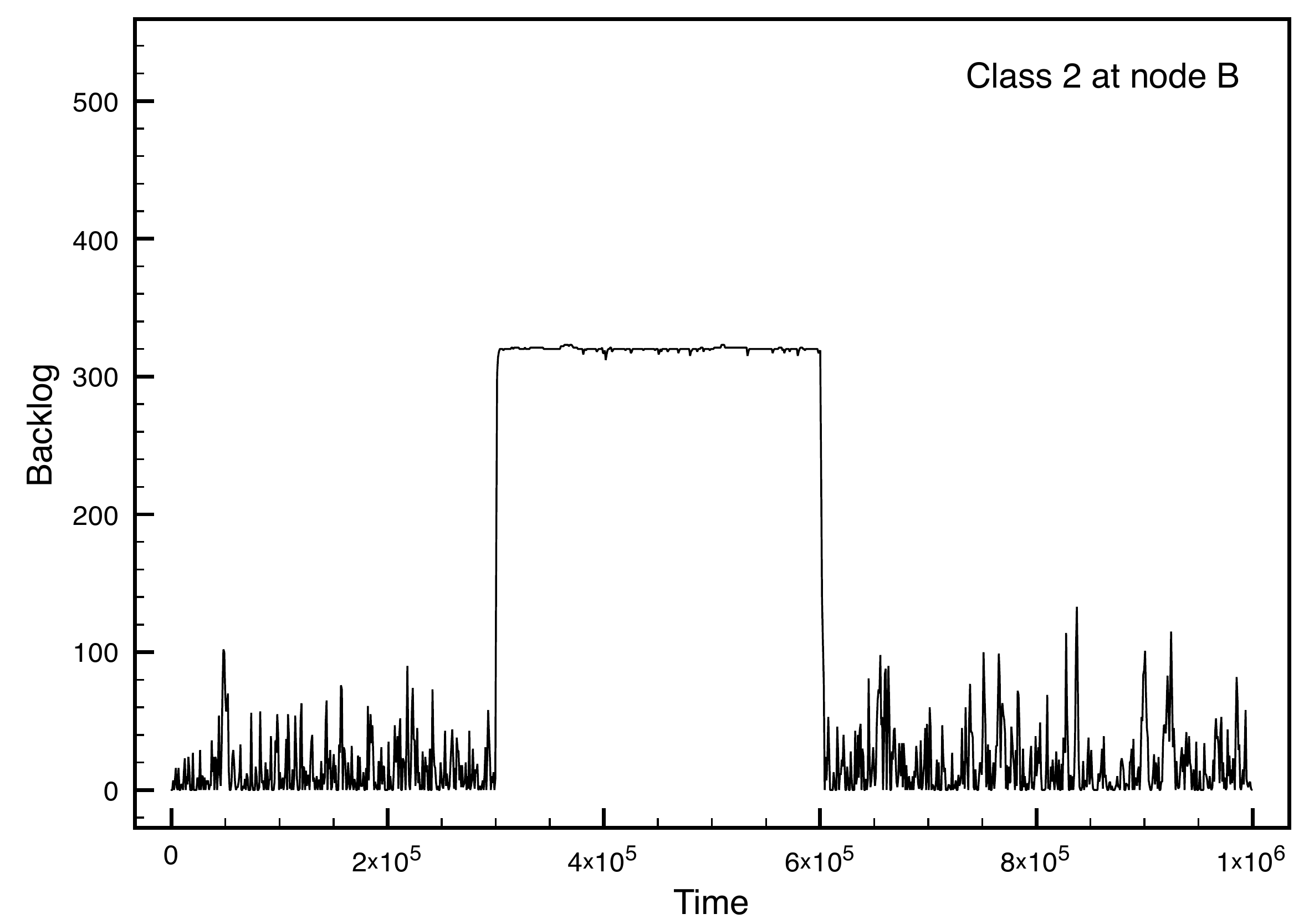}} \\
\subfigure[Queue $Q_{A}^{(2)}(t)$]{\includegraphics[width=1.7in]{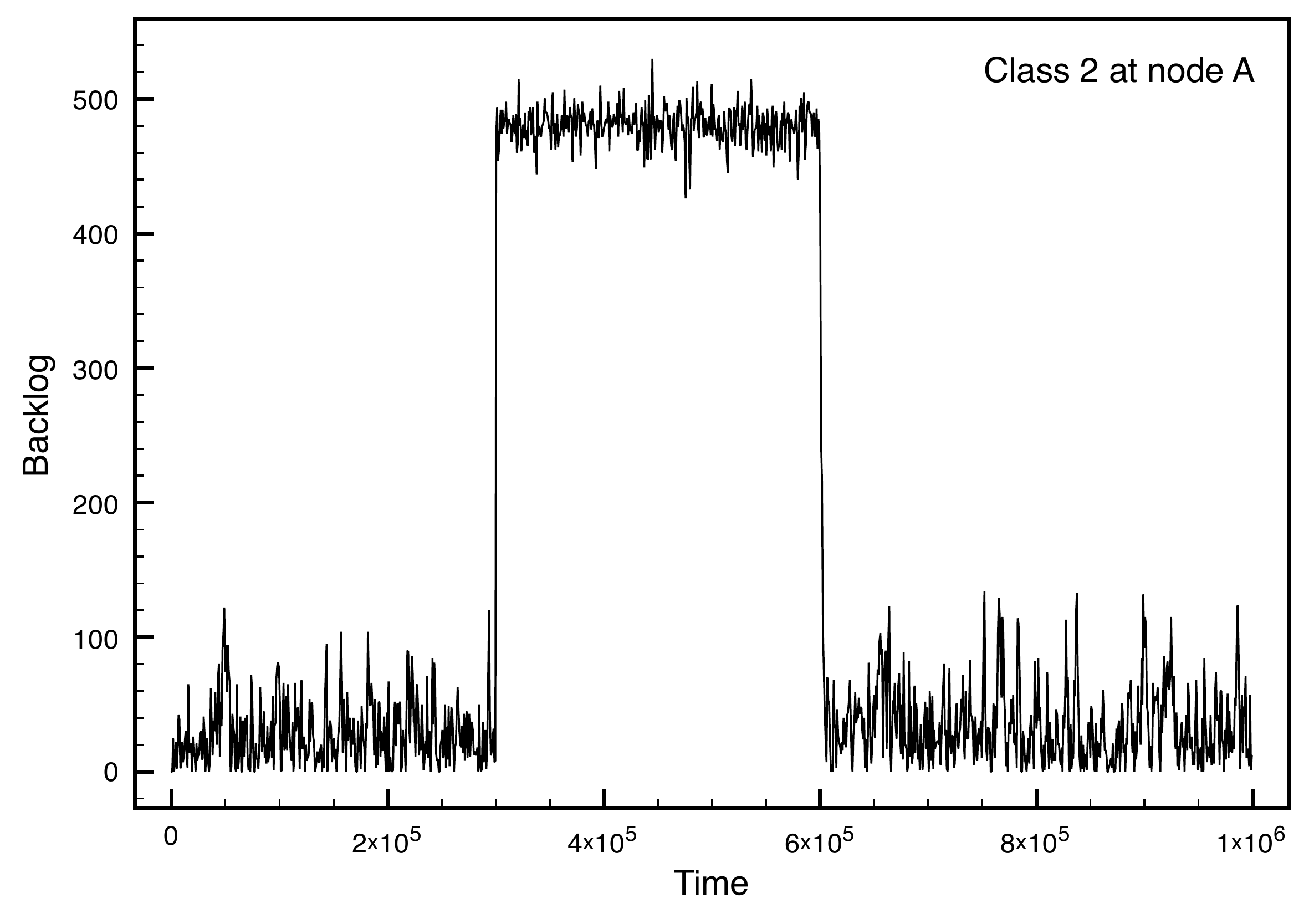}}
\subfigure[Queue $Q_{A}^{(3)}(t)$]{\includegraphics[width=1.7in]{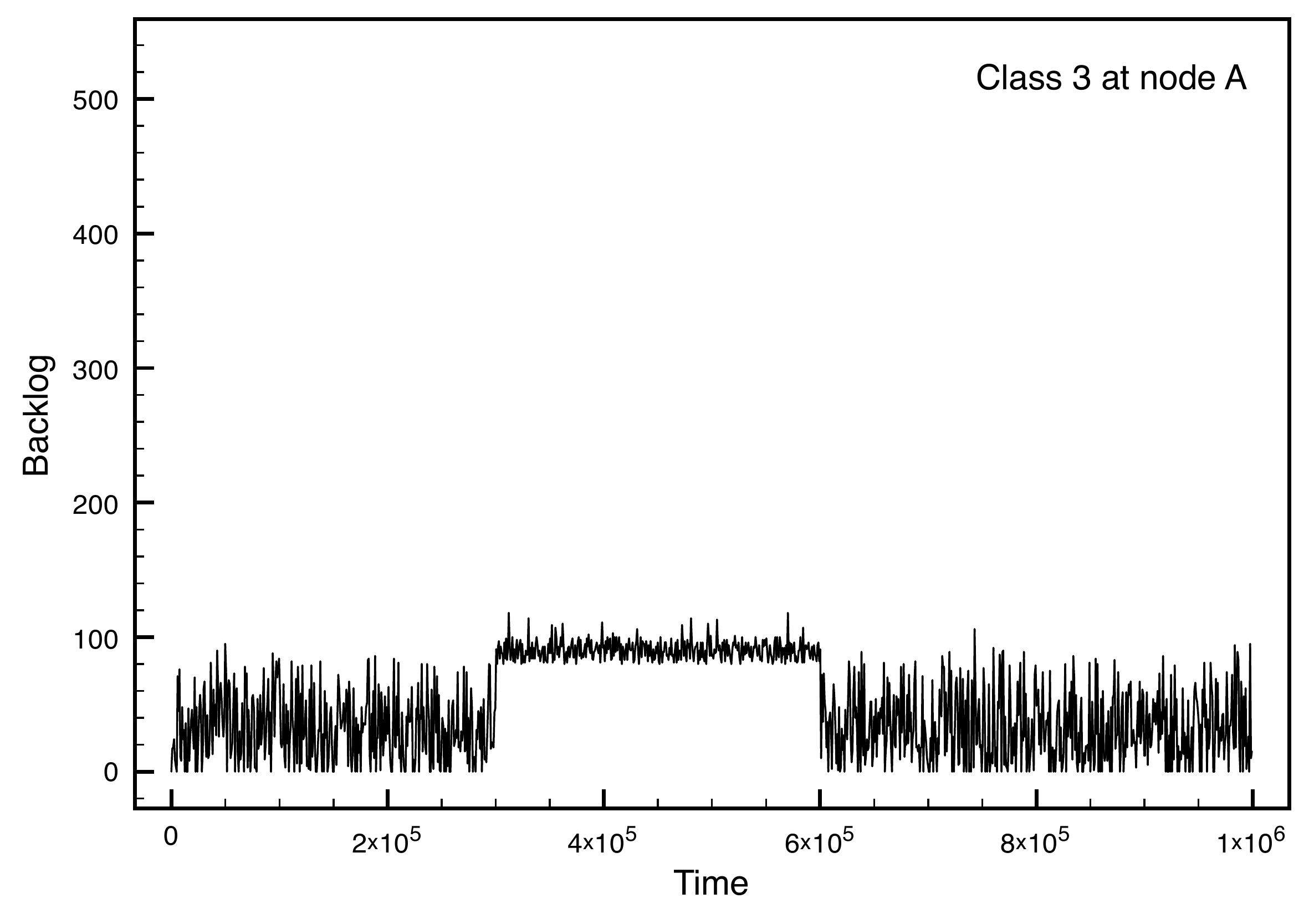}}
\caption{The queue processes under the $\ORA$ policy with time-varying arrival rates that temporarily overload the network.}
\label{fig:103}
\end{figure}

We use the following parameters here: $V=100$, $\amax=20$, $\dmax= \amax + \mumaxin = 21$, and $(\theta_{1}, \theta_{2}, \theta_{3}) = (3, 5, 1)$. Table~\ref{table:603} shows the near-optimal throughput performance of the $\ORA$ policy. Figure~\ref{fig:103} shows the sample paths of the queue processes $Q_{B}^{(1)}(t)$, $Q_{B}^{(2)}(t)$, $Q_{A}^{(2)}(t)$, and $Q_{A}^{(3)}(t)$ in the simulation. Clearly the queues suddenly build up when the network enters the overload interval $\mathcal{T}$, but the backlogs are kept close to the upper bound $\Qmax^{(c)} = V\theta_{c} + 2\dmax$ without blowing up.
\begin{table}
\caption{The throughput performance of the $\ORA$ policy under time-varying arrival rates.}
\label{table:603}
\centering
\small
\begin{tabular}{lll}
\toprule
time interval & throughput & optimal \\ 
& in this interval & value \\ \otoprule
$[0, 3\cdot10^{5})$ & $(.797, .097, .771)$ & $(.8, .1, .8)$ \\
$[3\cdot10^{5}, 6\cdot10^{5})$ & $(.001, .998, 0)$ & $(0, 1, 0)$ \\
$[6\cdot10^{5}, 10^{6})$ & $(.798, .102, .772)$ & $(.8, .1, .8)$ \\
\bottomrule
\end{tabular}
\end{table}

\section{Utility-Optimal Control} \label{sec:501}

We solve the general utility maximization problem~\eqref{eq:207}-\eqref{eq:208} with a network control policy very similar to the $\ORA$ policy in the previous section except for an additional flow control mechanism.

\subsection{Virtue Queue} \label{sec:1202}

In Section~\ref{sec:603} we formulate the equivalent optimization problem~\eqref{eq:901}-\eqref{eq:902} that involves maximizing $\sum_{c\in\cC} h_{c}(\nu_{c})$ subject to $r_{c} = \nu_{c}$ for all classes $c$, where $\nu_{c}$ are auxiliary control variables and $r_{c}$ is the throughput of class $c$. To enforce the constraint $r_{c}=\nu_{c}$, we construct a virtual queue $Z_{c}(t)$ in which $r_{c}$ is the virtual arrival rate and $\nu_{c}$ is the time-average service rate. Let $\tmu_{ad_{c}}^{(c)}(t)$ be the class $c$ packets received by node $d_{c}$ over link $(a, d_{c})$; we have
\[
\tmu_{ad_{c}}^{(c)}(t) = \min\left[Q_{a}^{(c)}(t), \mu_{ad_{c}}^{(c)}(t)\right].
\]
The arrivals to the virtual queue $Z_{c}(t)$ in a slot are the total  class $c$ packets delivered in that slot, namely, $\sum_{a} \tmu_{ad_{c}}^{(c)}(t)$. Let $\nu_{c}(t)$ be the allocated virtual service rate at $Z_{c}(t)$ in slot $t$. The virtual queue $\Zc(t)$ is located at the receiver $d_{c}$ and evolves  according to
\begin{equation} \label{eq:215}
\Zc(t+1) = \big[ \Zc(t) - \nuc(t) \big]^{+}+ \sum_{a} \tmu_{ad_{c}}^{(c)}(t).
\end{equation}
Assume initially $Z_{c}(0)=0$ for all classes $c$.   It is well known that if queue $Z_{c}(t)$ is stable then $r_{c} \leq \nu_{c}$. But we are interested in the stronger relationship that stabilizing $Z_{c}(t)$ leads to $r_{c} = \nu_{c}$. To make it happen, it suffices to guarantee that queue $Z_{c}(t)$ wastes as few service opportunities as possible, so that the time-average allocated service rate $\nu_{c}$ is approximately equal to the throughput out of queue $Z_{c}(t)$. For this, we need two conditions:
\begin{enumerate}
\item The queues $Z_{c}(t)$ usually have more than enough (virtual) data to serve.
\item When $Z_{c}(t)$ does not have sufficient data, the allocated service rate $\nu_{c}(t)$ is made arbitrarily small.
\end{enumerate}
To attain the first condition, we use an exponential-type Lyapunov function that bounds the virtual backlog process~$Z_{c}(t)$ away from zero (and centers it around a parameter $Q>0$). The second condition is attained by a proper choice of the parameters $\theta_{c}$ to be decided later.

\subsection{Control Policy}

The following policy, constructed in Appendix~\ref{appendix:501} and~\ref{appendix:202}, solves the general utility maximization problem~\eqref{eq:207}-\eqref{eq:208}.

\noindent \rule[0.05in]{3.5in}{0.01in}

Utility-Optimal Overload-Resilient Algorithm ($\UORA$)

\noindent \rule[0.05in]{3.5in}{0.01in}

\emph{Parameter Selection:} Choose positive parameters $\numax$, $w$,  $V$, $Q$, and $\{\theta_{c}, c\in \cC\}$ to be discussed shortly. Assume initially $\Qnc(0) = Z_{c}(0) =0$ and $\Dnc(0)=V \theta_{c}$.

\emph{Packet Dropping:} Same as the $\ORA$ policy.

\emph{Backpressure Routing:} Same as the $\ORA$ policy, except that the differential backlog over each link $l = (a, d_{c}) \in \cL$ connected to a receiver $d_{c}$ is modified as:
\begin{equation} \label{eq:1901}
W_{l}^{(c)}(t) = Q_{a}^{(c)}(t) - Q_{d_{c}}^{(c)}(t),
\end{equation}
where we abuse the notation by redefining
\begin{equation} \label{eq:610}
Q_{d_{c}}^{(c)}(t) = \begin{cases} w\,e^{w(Z_{c}(t)-Q)} & \text{if $Z_{c}(t)\geq Q$} \\ -w\,e^{w(Q-Z_{c}(t))} & \text{if $Z_{c}(t) < Q$} \end{cases}
\end{equation}
for all classes $c$. The exponential form of $Q_{d_{c}}^{(c)}(t)$ is a result of using exponential-type Lyapunov functions. We emphasize that here $Q_{d_{c}}^{(c)}(t)$ has nothing to do with real data buffered at the receivers (which must be zero); it is just a function of the virtual queue backlog $Z_{c}(t)$ that gives us the ``desired force'' in the form of differential backlog in~\eqref{eq:1901} to pull or push-back data in the network. Thus, unlike standard back-pressure routing that has $Q_{d_{c}}^{(c)}(t)=0$, here we use $Q_{d_{c}}^{(c)}(t)$ as part of the receiver-based flow control mechanism.

\emph{Receiver-Based Flow Control:} At a destination $d_{c}$, choose the virtual service rate $\nu_{c}(t)$ of queue $Z_{c}(t)$ as the solution to
\begin{align}
\text{maximize} &\quad V h_{c}\big(\nu_{c}(t)\big) + \nu_{c}(t) \, Q_{d_{c}}^{(c)}(t) \label{eq:406} \\
\text{subject to} &\quad 0\leq \nu_{c}(t)\leq \numax \label{eq:407}
\end{align}
where $h_{c}(x) = g_{c}(x) - \theta_{c}\, x$.

\emph{Queue Update:} Update queues $\Qnc(t)$, $\Dnc(t)$, and $Z_{c}(t)$ according to~\eqref{eq:201},~\eqref{eq:229}, and~\eqref{eq:215}, respectively, in every slot.

\noindent \rule[0.05in]{3.5in}{0.01in}

\subsection{Choice of Parameters} \label{sec:1201}

We briefly discuss how the parameters in the $\UORA$ policy are chosen. Let $\epsilon >0$ be a small constant which affects the performance of the $\UORA$ policy (cf.~\eqref{eq:906}). In~\eqref{eq:406}-\eqref{eq:407}, we need the parameter $\numax$ to satisfy $\numax \geq \max_{c\in\cC} r_{c}^{*} + \epsilon/2$, where $ (r_{c}^{*})$ is solution to the utility maximization~\eqref{eq:207}-\eqref{eq:208} (one feasible choice of $\numax$ is the sum of capacities of all links connected to the receivers plus $\epsilon/2$). This choice of $\numax$ ensures that queue $Z_{c}(t)$ can be stabilized when its virtual arrival rate is the optimal throughput $r_{c}^{*}$. Due to technical reasons, we define
$\deltamax \triangleq \max[\numax, \mumaxin]$ and choose the parameter
\[
w \triangleq \frac{\epsilon}{\deltamax^{2}} e^{-\epsilon/\deltamax}
\]
in~\eqref{eq:610}. The parameter $Q$ (see~\eqref{eq:610}) is used to bound the queues $Z_{c}(t)$ away from zero and center them around $Q$; for technical reasons, we need $Q \geq \numax$. The parameters $\theta_{c}$ are chosen to satisfy $h_{c}'(x) = g_{c}'(x) - \theta_{c} \leq 0$ for all $x\geq\epsilon$. This ensures that, when $Z_{c}(t) < Q$, its virtual service rate $\nu_{c}(t)$ as the solution to~\eqref{eq:406}-\eqref{eq:407} is less than or equal to $\epsilon$, attaining the second condition mentioned in Section~\ref{sec:1202} to equalize the arrival rate and the time-average service rate of the virtual queue $Z_{c}(t)$ (see Lemma~\ref{lem:403} in Appendix~\ref{appendix:502}). The parameter $V$ captures the tradeoff between utility and buffer sizes to be shown shortly and should be chosen large; for technical reasons, we need $V$ to satisfy $V\theta_{c} + 2\dmax \geq w$.

\subsection{Performance Analysis}

\begin{lem} \label{lem:502}
In the $\UORA$ policy, queues $\Qnc(t)$, $\Dnc(t)$, and $Z_{c}(t)$ are deterministically bounded by
\[
\Qnc(t)\leq \Qmax^{(c)},\quad \Dnc(t)\leq \Dmax^{(c)},\quad Z_{c}(t)\leq \Zmax^{(c)}
\]
for all $t$, $c$, and $n$, 
where $\Qmax^{(c)}$ and $\Dmax^{(c)}$ are defined in~\eqref{eq:503} and $\Zmax^{(c)}$ is defined as
\begin{equation} \label{eq:303}
\Zmax^{(c)} \triangleq Q + \frac{1}{w} \log \left(\frac{V \theta_{c} + 2\dmax}{w}\right) + \mumaxin.
\end{equation}
\end{lem}

\begin{IEEEproof}[Proof of Lemma~\ref{lem:502}]
See Appendix~\ref{appendix:602}.
\end{IEEEproof}

\begin{thm} \label{thm:502}
The $\UORA$ policy yields the limiting utility that satisfies
\begin{equation} \label{eq:906}
\sum_{c} g_{c}(\overline{r_{c}}) \geq \sum_{c} g_{c}(r_{c}^{*}) - \frac{B_{1}}{V} - \frac{3\epsilon}{2} \sum_{c} (m_{c} + \theta_{c}),
\end{equation}
where $\overline{r_{c}}$ is defined in~\eqref{eq:504}, $(r_{c}^{*})$ is the throughput vector that solves the utility maximization problem~\eqref{eq:207}-\eqref{eq:208}, and $B_{1}$ is a finite constant (defined in~\eqref{eq:1601}).
\end{thm}

\begin{IEEEproof}
See Appendix~\ref{appendix:502}.
\end{IEEEproof}

Theorem~\ref{thm:502} shows that the performance gap from the optimal utility can be made arbitrarily small by choosing a large $V$ and a small $\epsilon$. The performance tradeoff of choosing a large $V$ is again on the required finite buffer size $\Qmax^{(c)} = V\theta_{c} + 2\dmax$.

\subsection{Simulation of the $\UORA$ Policy}

We conduct two sets of simulations.

\subsubsection{On the $3$-node network in Fig.~\ref{fig:101}}
The goal is to provide proportional fairness to the three classes of traffic; equivalently we maximize the objective function $\log(r_{1})+\log(r_{2})+\log(r_{3})$. Each directed link $(A, B)$ and $(B, C)$ has the capacity of one packet per slot. The arrival process for each class is that, in every slot, $20$ packets arrive to the network with probability $0.1$ and zero packets arrive otherwise. The arrival rate vector is $(2, 2, 2)$, which overloads the network.  In this network setting,  due to symmetry the optimal throughput for class $1$ is equal to that of class $3$, which is the solution to the simple convex program
\[
\text{maximize:} \quad 2\log(x) + \log (1-x), \ \text{subject to:} \quad 0\leq x\leq 1.
\]
The optimal throughput vector is $(2/3, 1/3, 2/3)$ and the optimal utility is $-1.91$.  

As explained in Section~\ref{sec:1201}, we choose the parameters of the $\UORA$ policy as follows. Let $\epsilon = 0.1$. To satisfy $\theta_{c} \geq 1/x$ for all $x\geq \epsilon$, we choose $\theta_{c} = 1/\epsilon = 10$ for all classes $c$. The value of $\mumaxin$ in the 3-node network is one. The optimal throughput vector satisfies $\max_{c} r_{c}^{*} = 1$ and we choose $\numax = 3$ (any value of $\numax$ greater than $\max_{c} r_{c}^{*} +\epsilon/2 = 1.05$ works). By definition $\deltamax = \max[\numax, \mumaxin] = 3$. In the arrival processes we have $\amax=20$. By Assumption~\ref{assu:601} we choose $\dmax = \amax + \mumaxin = 21$. Let $Q=1000$.

We simulate the $\UORA$ policy for different values of $V$. The simulation time is $10^{6}$ slots. The near-optimal throughput performance is given in Table~\ref{table:605}. Table~\ref{table:606} shows the maximum backlog in each queue $\Qnc(t)$ during the simulation. Consistent with Lemma~\ref{lem:502}, the maximum backlog is bounded by $\Qmax^{(c)} = V\theta_{c} + 2\dmax = 10V+42$.
\begin{table}
\caption{The throughput performance of the $\UORA$ policy in the $3$-node network}
\label{table:605}
\centering
\small
\begin{tabular}{rcccc}
\toprule
$V$ & $r_{1}$ & $r_{2}$ & $r_{3}$ & $\sum \log (r_{c})$ \\ \otoprule
$10$ & $.522$ & $.478$ & $.522$ & $-2.038$ \\
$20$ & $.585$ & $.415$ & $.585$ & $-1.952$ \\
$50$ & $.631$ & $.369$ & $.631$ & $-1.918$ \\
$100$ & $.648$ & $.352$ & $.647$ & $-1.912$ \\ \otoprule
optimal & $.667$ & $.333$ & $.667$ & $-1.910$ \\
\bottomrule
\end{tabular}
\end{table}
\begin{table}
\caption{Maximum backlog in queues under the $\UORA$ policy in the $3$-node network}
\label{table:606}
\centering
\small
\begin{tabular}{rccccc}
\toprule
$V$ & $Q_{B}^{(1)}(t)$ & $Q_{B}^{(2)}(t)$ & $Q_{A}^{(2)}(t)$ & $Q_{A}^{(3)}(t)$ & $\Qmax^{(c)}$ \\ \otoprule
$10$ & $140$ & $97$ & $137$ & $137$ & $142$ \\
$20$ & $237$ & $187$ & $240$ & $236$ & $242$ \\
$50$ & $539$ & $441$ & $538$ & $540$ & $542$ \\
$100$ & $1036$ & $865$ & $1039$ & $1039$ & $1042$\\
\bottomrule
\end{tabular}
\end{table}

\subsubsection{On the tree network in Fig.~\ref{fig:102}}
Consider providing max-min fairness to the three classes of traffic in Fig.~\ref{fig:102}. Each link has the capacity of one packet per slot. Each one of the four arrival processes has 20 packets arriving in a slot with probability $0.1$ and zero packets arrive otherwise. The arrival rates are $(2, 2, 2, 2)$, which overloads the network. The optimal throughput for the three classes is easily seen to be $(2/3, 2/3, 2/3)$, where each flow of class $1$ contributes equally in that class.

We approximate max-min fairness by using the $\alpha$-fairness functions $g_{c}(x) = x^{1-\alpha}/(1-\alpha)$ with a large value of $\alpha=100$. The utility maximization becomes:
\begin{align*}
\text{maximize} &\quad \frac{-1}{99} \big[(r_{1A}+r_{1C})^{-99} + (r_{2B})^{-99} +(r_{3D})^{-99} \big] \\
\text{subject to} &\quad \text{$(r_{1A}, r_{1C}, r_{2B}, r_{3D})$ feasible in Fig.~\ref{fig:102}},
\end{align*}
where $r_{1A}$ is the throughput of class $1$ flow originating from node $A$; the other variables are similarly defined. 

According to Section~\ref{sec:1201}, we choose the parameters of the $\UORA$ policy as follows. We require $\theta_{c} \geq x^{-100}$ for all $x\geq \epsilon$. For convenience, let us choose $\theta_{c} = \epsilon =1$ for all classes $c$. The optimal throughput vector satisfies $\max_{c} r_{c}^{*} = 2$, achieved when the network always serves class $1$. We choose $\numax = 4$ (any value of $\numax$ greater than $\max_{c} r_{c}^{*} + \epsilon/2 = 2.5$ works). We observe from Fig.~\ref{fig:102} that $\mumaxin =2$, and we have $\deltamax = \max[\numax,\mumaxin]=4$. We have $\amax=20$ in the arrival processes and by Assumption~\ref{assu:601} we choose $\dmax = \amax + \mumaxin = 22$. Let $Q=100$.

We simulate the $\UORA$ policy for different values of $V$ and each simulation takes $10^{6}$ slots. The near-optimal performance of the $\UORA$ policy is given in Table~\ref{table:607}.
\begin{table}[tp]
\caption{The throughput performance of the $\UORA$ policy in the tree network}
\label{table:607}
\centering
\small
\begin{tabular}{rccc}
\toprule
$V$ & $r_{1}$ & $r_{2}$ & $r_{3}$ \\ \otoprule
$10$ & $.200$ & $.100$ & $.100$  \\
$20$ & $.364$ & $.206$ & $.205$  \\
$30$ & $.661$ & $.650$ & $.651$  \\
$50$ & $.667$ & $.667$ & $.667$  \\ \otoprule
optimal & $.667$ & $.667$ & $.667$  \\
\bottomrule
\end{tabular}
\end{table}
 
\section{Conclusion}

In this paper we develop a receiver-based flow control and an in-network packet dropping strategy to cope with network overload.  Our scheme is robust to uncooperative users who do not employ source-based flow control and malicious users that intentionally overload the network.  A novel feature of our scheme is a receiver-based backpressure/push-back mechanism that regulates data flows at the granularity of traffic classes, where packets can be classified based on aggregates of data flows.  This is in contrast to source-based schemes that can only differentiate between source-destination pairs. We show that when the receiver-based flow control scheme is combined with a threshold-based packet dropping policy at internal network nodes, optimal utility can be achieved.





\appendices
\section{Proof of Lemma~\ref{lem:201}} \label{appendix:201}

First we show~\eqref{eq:202}-\eqref{eq:204} are necessary conditions. Given a control policy, let $F_{ab}^{(c)}(t)$ be the amount of class $c$ packets transmitted over link $(a,b)$ in the interval $[0, t]$, and $\Qnc(t)$ be the class $c$ packets queued at node $n$ at time $t$. From the fact that the difference between incoming and outgoing packets at a node in $[0, t]$ is equal to the queue backlog at time $t$, we have
\begin{equation} \label{eq:205}
\sum_{\tau=0}^{t-1} \Anc(\tau) + \sum_{a} F_{an}^{(c)}(t) = \Qnc(t) + \sum_{b} F_{nb}^{(c)}(t),
\end{equation}
which holds for all nodes $n\neq d_{c}$ for each class $c\in \cC$. Taking expectation and time average of~\eqref{eq:205}, we obtain
\begin{equation} \label{eq:1501}
\lambdanc + \frac{1}{t} \expect{\sum_{a} F_{an}^{(c)}(t)} = \frac{\expect{\Qnc(t)}}{t} + \frac{1}{t} \expect{\sum_{b} F_{nb}^{(c)}(t)}.
\end{equation}
The link capacity constraints lead to
\begin{equation} \label{eq:206}
\frac{1}{t} \sum_{c} \expect{F_{ab}^{(c)}(t)} \leq \mumax^{ab}, \quad\forall\, (a,b)\in \cL.
\end{equation}
Consider the sequences $\mathbb{E}\big[F_{ab}^{(c)}(t)\big]/t$ and $\mathbb{E}\big[\Qnc(t)\big]/t$ indexed by $t$ in~\eqref{eq:1501}. For each $c\in\cC$ and $(a,b)\in\cL$, the sequence $\{\mathbb{E}\big[F_{ab}^{(c)}(t)\big]/t, t\in\Z^{+}\}$ is bounded because the capacity of each link is bounded. It follows from~\eqref{eq:1501} that the sequence $\{\mathbb{E}\big[\Qnc(t)\big]/t, t\in\Z^{+}\}$ is also bounded. There is a subsequence $\{t_{k}\}$ such that limit points $f_{ab}^{(c)}$ and $\qnc$ exist and satisfy, as $k\to\infty$,
\begin{align}
\frac{1}{t_{k}}  \expect{F_{ab}^{(c)}(t_{k})} &\to f_{ab}^{(c)}, \quad \forall\,c\in\cC,\ \forall\, (a,b)\in\cL,  \label{eq:601} \\
\frac{1}{t_{k}}\expect{\Qnc(t_{k})} &\to \qnc, \quad \forall\,c\in\cC,\ \forall\,n\neq d_{c}. \label{eq:602}
\end{align}
Applying~\eqref{eq:601}-\eqref{eq:602} to~\eqref{eq:1501}-\eqref{eq:206} results in~\eqref{eq:202} and~\eqref{eq:203}. Define the throughput $r_{c}$ of class $c$ as
\begin{equation} \label{eq:603}
r_{c} \triangleq \liminf_{t\to\infty} \frac{1}{t} \sum_{a: (a,d_{c})\in\cL} \expect{F_{ad_{c}}^{(c)}(t)}.
\end{equation}
The inequality in~\eqref{eq:204} follows~\eqref{eq:603} and~\eqref{eq:601}. The equality in~\eqref{eq:204} results from summing~\eqref{eq:202} over $n\neq d_{c}$.

To show the converse, it suffices to show that every interior point of $\Lambda_{\blambda}$ is achievable. Let $(r_{c})_{c\in\cC}$ be an interior point of $\Lambda_{\blambda}$, i.e., there exists $\epsilon \in (0, 1)$ such that $(r_{c} + \epsilon)_{c\in\cC} \in\Lambda_{\blambda}$.  There exist corresponding flow variables $f_{ab}^{(c)}$ and $\qnc$ such that
\begin{gather}
\lambdanc + \sum_{a} f_{an}^{(c)} = q_{n}^{(c)} + \sum_{b} f_{nb}^{(c)}, \ \sum_{c} f_{ab}^{(c)} \leq \mumax^{ab}, \label{eq:1502} \\
r_{c} +\epsilon \leq \sum_{n} \lambdanc - \sum_{n} \qnc. \label{eq:1503}
\end{gather}
In the flow system~\eqref{eq:1502}-\eqref{eq:1503}, by removing the subflows that contribute to queue overflows, we obtain new reduced flow variables $\hat{\lambda}_{n}^{(c)}$  and $\hat{f}_{ab}^{(c)}$ such that $0\leq\hat{\lambda}_{n}^{(c)} \leq \lambdanc$, $0\leq\hat{f}_{ab}^{(c)} \leq f_{ab}^{(c)}$, and
\begin{gather}
\hat{\lambda}_{n}^{(c)} + \sum_{a} \hat{f}_{an}^{(c)} = \sum_{b} \hat{f}_{nb}^{(c)}, \quad \forall\, c,\ \forall n\neq d_{c}, \label{eq:904} \\
\quad \sum_{c} \hat{f}_{ab}^{(c)} \leq \mumax^{ab}, \quad \forall (a,b)\in\cL, \label{eq:604} \\
r_{c} + \epsilon \leq \sum_{n} \hat{\lambda}_{n}^{(c)}, \quad \forall c\in\cC. \label{eq:903}
\end{gather}
Define
\[
r_{n}^{(c)} \triangleq \hat{\lambda}_{n}^{(c)} \left[1 - \frac{\epsilon/2}{\sum_{n} \hat{\lambda}_{n}^{(c)}} \right] \frac{r_{c}}{\sum_{n} \hat{\lambda}_{n}^{(c)}-\epsilon/2}.
\]
From~\eqref{eq:903}, we have $\sum_{n} \hat{\lambda}_{n}^{(c)} > \epsilon/2$ and $\sum_{n} \hat{\lambda}_{n}^{(c)} - \epsilon/2 \geq r_{c}$. It is not difficult to check that $r_{n}^{(c)} \geq 0$, $r_{n}^{(c)} < \hat{\lambda}_{n}^{(c)}$, and $\sum_{n} r_{n}^{(c)} = r_{c}$. Combined with~\eqref{eq:904}-\eqref{eq:604}, we obtain
\[
r_{n}^{(c)} + \sum_{a} \hat{f}_{an}^{(c)} < \sum_{b} \hat{f}_{nb}^{(c)}, \quad \sum_{c} \hat{f}_{ab}^{(c)} \leq \mumax^{ab},
\]
where the first inequality is a strict one. These inequalities show that the rate matrix $(r_{n}^{(c)})$ is an interior point of the network capacity region $\Lambda$ in Definition~\ref{defn:101}, and therefore is achievable by a control policy, such as the back-pressure policy~\cite{NMR05}. Therefore, the aggregate rate vector $(r_{c})$, where $r_{c} = \sum_{n} r_{n}^{(c)}$, is also achievable.

\section{Proof of Lemma~\ref{lem:402}} \label{appendix:601}

We prove Lemma~\ref{lem:402} by induction. First we show $\Dnc(t)$ is deterministically bounded. Assume $\Dnc(t) \leq\Dmax^{(c)}$ for some $t$, which holds at $t=0$ because we let $\Dnc(0) = V\theta_{c}$. Consider the two cases:
\begin{enumerate}
\item If $D_{n}^{(c)}(t) \leq V\theta_{c}$, then from~\eqref{eq:229} we obtain
\[
D_{n}^{(c)}(t+1) \leq D_{n}^{(c)}(t) + \widetilde{d}_{n}^{(c)}(t) \leq V\theta_{c}+ \dmax = \Dmax^{(c)},
\]
where the second inequality uses~\eqref{eq:1203}.
\item If $D_{n}^{(c)}(t) > V\theta_{c}$, then the $\ORA$ policy chooses $\varphi_{n}^{(c)}(t) = \dmax$ at queue $\Dnc(t)$ and we have
\begin{align*}
\Dnc(t+1) &\leq \big[\Dnc(t) - \dmax\big]^{+} + \dmax  \\
&\leq \max[\Dnc(t), \dmax] \leq \Dmax^{(c)},
\end{align*}
where the last inequality uses the induction assumption.
\end{enumerate}
From these two cases,  we obtain $D_{n}^{(c)}(t+1) \leq \Dmax^{(c)}$. 

Similarly, we show $\Qnc(t)$ is bounded. Assume $Q_{n}^{(c)}(t) \leq \Qmax^{(c)}$ for some $t$, which holds at $t=0$ because we let $\Qnc(0)=0$. Consider the two cases:
\begin{enumerate}
\item If $Q_{n}^{(c)}(t) \leq \Dmax^{(c)}$, then from~\eqref{eq:201} we get
\[
Q_{n}^{(c)}(t+1) \leq Q_{n}^{(c)}(t) + \amax + \mumaxin \leq \Dmax^{(c)} + \dmax = \Qmax^{(c)},
\]
where the second inequality uses Assumption~\ref{assu:601}.
\item If $Q_{n}^{(c)}(t) > \Dmax^{(c)} \geq D_{n}^{(c)}(t)$, the $\ORA$ policy chooses $d_{n}^{(c)}(t) = \dmax$ at queue $\Qnc(t)$ and we have
\[
\begin{split}
Q_{n}^{(c)}(t+1) &\leq Q_{n}^{(c)}(t) - \dmax + \amax + \mumaxin \\ &\leq Q_{n}^{(c)}(t) \leq \Qmax^{(c)},
\end{split}
\]
where the third inequality uses induction assumption.
\end{enumerate}
We conclude that $Q_{n}^{(c)}(t+1)  \leq \Qmax^{(c)}$. 

Finally, we show $\Dnc(t) \geq V\theta_{c} - \dmax$ for all slots. Assume this is true at some time $t$;  this holds when $t=0$ because we assume $\Dnc(0) = V\theta_{c}$. Consider the two cases:
\begin{enumerate}
\item If $\Dnc(t) \leq V\theta_{c}$, the $\ORA$ policy chooses $\varphinc(t) = 0$ at queue $\Dnc(t)$ and we have
\[
\Dnc(t+1) \geq \Dnc(t) \geq V\theta_{c} - \dmax
\]
by induction assumption.
\item If $\Dnc(t) > V\theta_{c}$, the $\ORA$ policy chooses $\varphinc(t) = \dmax$ and we have
\[
\Dnc(t+1) \geq \Dnc(t) - \dmax > V\theta_{c} - \dmax.
\]
\end{enumerate}
The proof is complete.

\section{} \label{appendix:501}

We construct a proper Lypuanov drift inequality that leads to the $\UORA$ policy. Let
\[
H(t) = \big(\Qnc(t); \Dnc(t); Z_{c}(t)\big)
\]
be the vector of all physical and virtual queues in the network. Using the parameters $w$ and $Q$ given in the policy, we define the Lyapunov function
\begin{multline*}
L\bigparen{H(t)} \triangleq \frac{1}{2} \sum_{c,n\neq d_{c}} \bigl[Q_{n}^{(c)}(t)\bigr]^{2} + \frac{1}{2}\sum_{c, n\neq d_{c}} \bigl[D_{n}^{(c)}(t)\bigr]^{2} \\
+ \sum_{c\in\cC} \Big(e^{w\bigparen{Z_{c}(t)-Q}} + e^{w\bigparen{Q-Z_{c}(t)}} \Big).
\end{multline*}
The last sum is a Lyapunov function whose value grows exponentially when $Z_{c}(t)$ deviates in both directions from the constant $Q$. This exponential-type Lyapunov function is useful for both stabilizing $Z_{c}(t)$ and guaranteeing  there is sufficient data in $Z_{c}(t)$. Such exponential-type Lyapunov functions are  previously used in~\cite{Nee06a} to study the optimal utility-delay tradeoff in wireless networks. We define the Lyapunov drift
\[
\Delta(t) \triangleq \expect{L\bigl(H(t+1)\bigr) - L\bigl(H(t)\bigr) | H(t)},
\]
where the expectation is with respect to all randomness in the system in slot $t$.

Define the indicator function $\onecR(t) = 1$ if $Z_{c}(t) \geq Q$ and $0$ otherwise; define $\onecL(t) = 1-\onecR(t)$. Define $\delta_{c}(t) = \nu_{c}(t) - \sum_{a: (a, d_{c})\in\cL} \mu_{ad_{c}}^{(c)}(t)$. The next lemma is proved in Appendix~\ref{appendix:202}.

\begin{lem} \label{lem:401}
The Lyapunov drift $\Delta(t)$ under any control policy satisfies
\begin{align}
&\Delta(t) -V \sum_{c} \expect{h_{c}\big(\nu_{c}(t)\big) \mid H(t)} \notag \\
 &+ V\sum_{nc}\theta_{c} \,\expect{\varphinc(t)\mid H(t)} \leq B + \sum_{nc} \Qnc(t) \lambdanc  \notag\\
&-\sum_{nc} \Qnc(t) \expect{\sum_{b} \mu_{nb}^{(c)}(t) + d_{n}^{(c)}(t) - \sum_{a} \mu_{an}^{(c)}(t) \mid H(t)} \notag \\
&- \sum_{nc} \Dnc(t) \expect{\varphinc(t)-\dnc(t) \mid H(t)} \notag \\
&-V \sum_{c} \expect{h_{c}\big(\nu_{c}(t)\big) \mid H(t)} + V\sum_{nc} \theta_{c}\, \expect{\varphinc(t)\mid H(t)} \notag \\
&- w \sum_{c}  \onecR(t) e^{w(Z_{c}(t)-Q)} \Big(\expect{\delta_{c}(t)\mid H(t)} - \frac{\epsilon}{2} \Big) \notag \\
&+ w \sum_{c} \onecL(t) e^{w(Q-Z_{c}(t))}\Big(\expect{ \delta_{c}(t) \mid H(t)} + \frac{\epsilon}{2} \Big),\label{eq:227}
\end{align}
where $B$ is a finite constant defined by
\[
\begin{split}
&B = \abs{\cN}\abs{\cC}\big[ (\mumaxout +\dmax)^{2}+(\amax+\mumaxin)^{2}\big] + \\
&\quad 2 \abs{\cN}\abs{\cC} \dmax^{2} +\abs{\cC}\big[ w(2\deltamax+\epsilon) + e^{w(\numax+\mumaxin)} + e^{wQ} \big].
\end{split}
\]
\end{lem}
By isolating decisions variables in~\eqref{eq:227}, it is not difficult to verify that the $\UORA$ policy observes the current network state $H(t)$ and minimizes the right-hand side of~\eqref{eq:227} in every slot.

\section{} \label{appendix:202}

We establish the Lyapunov drift inequality in~\eqref{eq:227}. Applying to~\eqref{eq:201} the fact that $[ (a-b)^{+}-c]^{+} = (a-b-c)^{+}$ for nonnegative reals $a$, $b$, and $c$, and using Lemma~\ref{lem:101} in Appendix~\ref{appendix:203}, we have
\begin{equation} \label{eq:126}
\begin{split}
&\frac{1}{2} \Big( \big[ Q_{n}^{(c)}(t+1)\big]^{2} - \big[Q_{n}^{(c)}(t)\big]^{2} \Big) \leq B_{Q} - \\
&Q_{n}^{(c)}(t)\, \Big( \sum_{b} \mu_{nb}^{(c)}(t) + d_{n}^{(c)}(t)  - A_{n}^{(c)}(t) - \sum_{a} \mu_{an}^{(c)}(t) \Big),
\end{split}
\end{equation}
where $B_{Q} = (\mumaxout + \dmax)^{2} + (\amax + \mumaxin)^{2}$ is a finite constant. From~\eqref{eq:229} and~\eqref{eq:1203}, we get
\[
\Dnc(t+1) \leq \big[ \Dnc(t) - \varphinc(t) \big]^{+}+ \dnc(t).
\]
Similarly, we obtain for queue $\Dnc(t)$
\begin{multline} \label{eq:127}
\frac{1}{2} \Big( \big[ D_{n}^{(c)}(t+1)\big]^{2} - \big[D_{n}^{(c)}(t)\big]^{2} \Big) \\
\leq B_{D} - D_{n}^{(c)}(t)\, \big( \varphi_{n}^{(c)}(t) - d_{n}^{(c)}(t) \big),
\end{multline}
where $B_{D} = 2\dmax^{2}$ is a finite constant.

\begin{lem} \label{lem:202}
Given fixed constants $\epsilon >0$ and $Q>0$, we define
\begin{equation} \label{eq:216}
w \triangleq \frac{\epsilon}{\deltamax^{2}} e^{-\epsilon/\deltamax}.
\end{equation}
Then
\begin{multline} \label{eq:222}
e^{w(Z_{c}(t+1)-Q)} - e^{w(Z_{c}(t)-Q)} \\ \leq e^{w(\numax + \mumaxin)} - w e^{w(Z_{c}(t)-Q)} \Big[\delta_{c}(t) - \frac{\epsilon}{2} \Big],
\end{multline}
\begin{multline} \label{eq:223}
e^{w (Q-Z_{c}(t+1))} - e^{w(Q-Z_{c}(t))}\\
\leq e^{wQ} + w e^{w (Q-Z_{c}(t))} \Big[\delta_{c}(t) + \frac{\epsilon}{2}\Big].
\end{multline}
\end{lem}
\begin{IEEEproof}
Define $\delta_{c}(t) = \nu_{c}(t) - \sum_{a} \mu_{ad_{c}}^{(c)}(t)$ and $\deltamax = \max(\numax, \mumaxin)$; we have $\abs{\delta_{c}(t)} \leq \deltamax$ and $\nu_{c}(t) \leq \deltamax$. Equation~\eqref{eq:215} yields
\[
Z_{c}(t+1) \leq \begin{cases} Z_{c}(t) - \delta_{c}(t) & \text{if $Z_{c}(t) \geq \numax$,} \\ \numax + \mumaxin & \text{if $Z_{c}(t) < \numax$.} \end{cases}
\]
Since $w>0$, we have
\begin{equation} \label{eq:1205}
e^{wZ_{c}(t+1)} \leq e^{w(\numax + \mumaxin)} + e^{wZ_{c}(t)} e^{-w\delta_{c}(t)},
\end{equation}
because the first term is bounded by the second term if $Z_{c}(t) < \numax$, and is bounded by the third term otherwise. Multiplying both sides by $\exp(-wQ)$ yields
\begin{multline} \label{eq:217}
e^{w(Z_{c}(t+1)-Q)} \leq e^{w(\numax + \mumaxin)} \\
+ e^{w(Z_{c}(t)-Q)} \Big[ 1 - w\delta_{c}(t) + \frac{w^{2}\deltamax^{2}}{2} e^{w\deltamax} \Big],
\end{multline}
where the last term follows the Taylor expansion of $e^{-w \delta_{c}(t)}$. If we have
\begin{equation} \label{eq:301}
\frac{w^{2}\deltamax^{2}}{2} e^{w\deltamax} \leq \frac{w\epsilon}{2},
\end{equation}
then plugging~\eqref{eq:301} into~\eqref{eq:217} leads to~\eqref{eq:222}. Indeed, by definition of $w$ in~\eqref{eq:216} we get
\[
\frac{w^{2}\deltamax^{2}}{2} e^{w\deltamax} = \frac{w\epsilon}{2} e^{w\deltamax-\epsilon/\deltamax} \leq \frac{w\epsilon}{2},
\]
which uses $w \leq \epsilon/\deltamax^{2} \Rightarrow \exp(w\deltamax) \leq \exp(\epsilon/\deltamax)$.

Also,~\eqref{eq:215} leads to
\[
Z_{c}(t+1) \geq Z_{c}(t) - \nu_{c}(t) + \sum_{a} \tmu_{ad_{c}}^{(c)}(t).
\]
Define the event $E$ as
\[
E: \text{if $Q_{n}^{(c)}(t) \geq \mumax^{nd_{c}}\ \forall\, n:(n,d_{c})\in \cL$,}
\]
i.e., all upstream nodes of $d_{c}$ have sufficient data to transmit, in which case we must have $\tmu_{ad_{c}}^{(c)}(t) = \mu_{ad_{c}}^{(c)}(t)$. It follows
\begin{equation} \label{eq:609}
Z_{c}(t+1) \geq \begin{cases} Z_{c}(t)-\delta_{c}(t) & \text{if event $E$ happens} \\ 0 & \text{otherwise}\end{cases}
\end{equation}
where the second case follows that queue $Z_{c}(t)$ is always non-negative. Similar as~\eqref{eq:1205}, from~\eqref{eq:609} we obtain
\begin{equation} \label{eq:302}
\begin{split}
&e^{-w Z_{c}(t+1)} \leq 1 + e^{-w Z_{c}(t)} e^{w\delta_{c}(t)} \\
&\leq 1+ e^{-w Z_{c}(t)} \Big[1+ w \delta_{c}(t) + \frac{w^{2}\deltamax^{2}}{2} e^{w\deltamax}\Big],
\end{split}
\end{equation}
in which we use the Taylor expansion of $e^{w\delta_{c}(t)}$. Bounding the last term of~\eqref{eq:302} by~\eqref{eq:301} and multiplying the inequality by  $\exp(wQ)$, we obtain~\eqref{eq:223}.
\end{IEEEproof}

Now, define the indicator function $\onecR(t) = 1$ if $Z_{c}(t) \geq Q$ and $0$ otherwise. Define $\onecL(t) = 1-\onecR(t)$. Then for any real number $a$,
\begin{equation} \label{eq:220}
\begin{split}
a\, e^{w(Z_{c}(t)-Q)} &= a (\onecR(t)+\onecL(t)) e^{w(Z_{c}(t)-Q)} \\
&\leq \abs{a} + a\,\onecR(t)\, e^{w(Z_{c}(t)-Q)}.
\end{split}
\end{equation}
Similarly, we have
\begin{equation} \label{eq:221}
\begin{split}
a\, e^{w(Q-Z_{c}(t))} \leq \abs{a} + a\,\onecL(t)\, e^{w(Q-Z_{c}(t))}.
\end{split}
\end{equation}
Using $a = -w (\delta_{c}(t)-\epsilon/2)$ in~\eqref{eq:220} and $a = w (\delta_{c}(t)+ \epsilon/2)$ in~\eqref{eq:221},  inequalities~\eqref{eq:222} and~\eqref{eq:223} lead to
\begin{multline} \label{eq:224}
e^{w(Z_{c}(t+1)-Q)} - e^{w(Z_{c}(t)-Q)} \leq e^{w(\numax + \mumaxin)} \\ + w(\deltamax+\frac{\epsilon}{2}) - w \onecR(t) e^{w(Z_{c}(t)-Q)} \Big[\delta_{c}(t) - \frac{\epsilon}{2} \Big],
\end{multline}
\begin{multline} \label{eq:225}
e^{w(Q-Z_{c}(t+1))} - e^{w(Q-Z_{c}(t))} \leq e^{wQ} + w(\deltamax+\frac{\epsilon}{2}) \\
+ w \onecL(t) e^{w(Q-Z_{c}(t))} \Big[\delta_{c}(t) + \frac{\epsilon}{2} \Big].
\end{multline}

Finally, summing~\eqref{eq:126} and~\eqref{eq:127} over $n\neq d_{c}$ for each $c\in \cC$, summing~\eqref{eq:224} and~\eqref{eq:225} over $c\in \cC$, and taking conditional expectation, we have
\begin{multline} \label{eq:226}
\Delta(t) \leq B + \sum_{nc} \Qnc(t) \lambdanc  \\
-\sum_{nc} \Qnc(t) \expect{\sum_{b} \mu_{nb}^{(c)}(t) + d_{n}^{(c)}(t) - \sum_{a} \mu_{an}^{(c)}(t) \mid H(t)} \\
- \sum_{nc} \Dnc(t) \expect{\varphinc(t)-\dnc(t) \mid H(t)} \\
- w \sum_{c} \onecR(t) e^{w(Z_{c}(t)-Q)} \Big(\expect{\delta_{c}(t)\mid H(t)} - \frac{\epsilon}{2} \Big) \\
+ w \sum_{c} \onecL(t) e^{w(Q-Z_{c}(t))}\Big(\expect{\delta_{c}(t)\mid H(t)} + \frac{\epsilon}{2} \Big),
\end{multline}
where $B$ is a finite constant defined as
\begin{align}
B &= \abs{\cN}\abs{\cC}(B_{Q}+B_{D}) \notag \\
&\quad + \abs{\cC} \Big[w (2\deltamax + \epsilon)  + e^{w(\numax + \mumaxin)} + e^{wQ}\Big]. \label{eq:1208}
\end{align}
The constants $B_{Q}$ and $B_{D}$ are defined right after~\eqref{eq:126} and~\eqref{eq:127}, respectively. Adding to both sides of~\eqref{eq:226}
\[
-V \sum_{c} \expect{h_{c}\big(\nu_{c}(t)\big) \mid H(t)} +  V\sum_{nc} \theta_{c}\, \expect{\varphinc(t)\mid H(t)}
\]
yields~\eqref{eq:227}.

\section{Proof of Lemma~\ref{lem:502}} \label{appendix:602}

The boundedness of the queues $\Qnc(t)$ and $\Dnc(t)$ follows the proof of Lemma~\ref{lem:402}. For the queues $Z_{c}(t)$, we again prove by induction. Suppose $Z_{c}(t) \leq \Zmax^{(c)}$, which holds at $t=0$ because we let $Z_{c}(0)=0$. Consider the two cases:
\begin{enumerate}
\item If $Z_{c}(t) \leq \Zmax - \mumaxin$, then by~\eqref{eq:215}
\[
Z_{c}(t+1) \leq Z_{c}(t) + \mumaxin \leq \Zmax^{(c)}.
\]
\item If $Z_{c}(t) >  \Zmax^{(c)}- \mumaxin$, then from~\eqref{eq:303} we obtain
\begin{align}
Z_{c}(t)-Q &> \Zmax^{(c)}- \mumaxin -Q \notag \\
&\stackrel{(a)}{=} \frac{1}{w} \log \left(\frac{V \theta_{c} + 2\dmax}{w}\right) \geq 0, \label{eq:1902}
\end{align}
where the last inequality follows our choice of $V$ satisfying $V\theta_{c} + 2\dmax \geq w$ as mentioned in Section~\ref{sec:1201}.
As a result, in~\eqref{eq:610} we have
\[
Q_{d_{c}}^{(c)}(t) = w e^{w(Z_{c}(t)-Q)} >V\theta_{c} + 2 \dmax = \Qmax^{(c)},
\]
where the inequality follows $(a)$ in~\eqref{eq:1902}. Since all queues $\Qnc(t)$ for $n\neq d_{c}$ are deterministically bounded by $\Qmax^{(c)}$, we have $Q_{d_{c}}^{(c)}(t) > \Qnc(t)$ for all nodes $n$ such that $(n,d_{c})\in \cL$. Consequently, the $\UORA$ policy does not transmit any class $c$ packets over the links $(n,d_{c})$ and the virtual queue $Z_{c}(t)$ has no arrivals. Therefore $Z_{c}(t+1) \leq Z_{c}(t) \leq \Zmax^{(c)}$.
\end{enumerate}
We conclude that $Z_{c}(t+1) \leq \Zmax^{(c)}$. The proof is complete.

\section{Proof of Theorem~\ref{thm:502}} \label{appendix:502}

Let the throughput vector $(r_{c}^{*})$ and the flow variables $(q_{n}^{(c)*}; f_{ab}^{(c)*})$ be the  optimal solution to the utility maximization~\eqref{eq:212}-\eqref{eq:605}. We consider the stationary policy that observes the current network state $H(t)$ and chooses in every slot:
\begin{enumerate}
\item $\varphincstar(t) = \dncstar(t) = q_{n}^{(c)*}$.
\item $\mu_{ab}^{(c)*}(t) = f_{ab}^{(c)*}$.
\item $\nu_{c}^{*}(t) = r_{c}^{*} + \frac{\epsilon}{2}$ if $Z_{c}(t) \geq Q$, and $\nu_{c}^{*}(t) = r_{c}^{*}$ otherwise.
\end{enumerate}
The first part is a feasible allocation because the value of $q_{n}^{(c)*}$ is at most $\amax + \mumaxin$, which is less than or equal to $\dmax$ by Assumption~\ref{assu:601}. The second part is feasible because the link rates $f_{nm}^{(c)*}$ satisfy the link capacity constraints~\eqref{eq:203}. The third part is feasible because we assume $\numax \geq \max_{c\in\cC} r_{c}^{*} + \epsilon/2$.

Under this stationary policy, by using the equality in~\eqref{eq:204} and using~\eqref{eq:213}, we have $\sum_{a} \mu_{ad_{c}}^{(c)*}(t) = r_{c}^{*}$ in every slot. As a result,
\begin{align*}
\expect{\delta_{c}^{*}(t) \mid Z_{c}(t) \geq Q} &= \mathbb{E}\Big[\nu_{c}^{*}(t) - \sum_{a} \mu_{ad_{c}}^{(c)*}(t)  \mid Z_{c}(t) \geq Q\Big] \\
&= \expect{\nu_{c}^{*}(t) - r_{c}^{*} \mid Z_{c}(t) \geq Q} = \frac{\epsilon}{2}, \\
\expect{\delta_{c}^{*}(t)\mid Z_{c}(t) < Q} &= 0.
\end{align*}
The $\UORA$ policy minimizes the right-hand side of~\eqref{eq:227} in every slot. Using the decisions in the stationary policy, we can upper bound~\eqref{eq:227} under the $\UORA$ policy by
\begin{align} 
&\Delta(t) -V \sum_{c} \expect{h_{c}(\nu_{c}(t)) \mid H(t)} + V\sum_{nc} \theta_{c}\,\expect{\varphinc(t)\mid H(t)} \notag \\
&\leq B + \abs{\cC} \frac{w\epsilon}{2} e^{wQ}- V \sum_{c} h_{c}(\nu_{c}^{*}(t)) + V \sum_{nc} \theta_{c}\, q_{n}^{(c)*} \notag \\
&\ -\sum_{nc} \Qnc(t) \Big( \sum_{b} f_{nb}^{(c)*} + q_{n}^{(c)*} - \lambdanc - \sum_{a} f_{an}^{(c)*}\Big) \notag \\
&= B_{1} - V \sum_{c} h_{c}(\nu_{c}^{*}(t)) + V \sum_{nc}\theta_{c}\, q_{n}^{(c)*}, \label{eq:404}
\end{align}
where we define 
\begin{equation} \label{eq:1209}
B_{1} = B + \frac{1}{2}\abs{\cC}w\epsilon e^{wQ}
\end{equation}
and $B$ is given in~\eqref{eq:1208}.  The last equality in~\eqref{eq:404} is because the variables $f_{ab}^{(c)*}$ and $q_{n}^{(c)*}$ satisfy the flow conservation constraints~\eqref{eq:202}.

The utility functions $g_{c}(\cdot)$ are assumed to have bounded derivatives with $\abs{g_{c}'(x)} \leq m_{c}$ for all $x\geq 0$. Thus, the utility functions $h_{c}(\cdot)$ have bounded derivatives with $\abs{h_{c}'(x)} \leq m_{c}+\theta_{c}$ for all $x\geq 0$. Using $\nu_{c}^{*}(t) \in \{r_{c}^{*}, r_{c}^{*}+\epsilon/2\}$, we have
\begin{equation} \label{eq:1206}
h_{c}(\nu_{c}^{*}(t)) \geq h_{c}\big(r_{c}^{*}\big) - \frac{\epsilon}{2} (m_{c}+\theta_{c}).
\end{equation}
If $\nu_{c}^{*}(t) = r_{c}^{*}$ then it is trivial. If $\nu_{c}^{*}(t) = r_{c}^{*} + \epsilon/2$ then by mean-value theorem there exists $y\in(r_{c}^{*}, r_{c}^{*} + \epsilon/2)$ such that 
\[
\begin{split}
h_{c}(\nu_{c}^{*}(t)) &= h_{c}\big(r_{c}^{*}\big) + \big(\nu_{c}^{*}(t)-r_{c}^{*}\big)\, h_{c}'(y) \\
&\geq h_{c}\big(r_{c}^{*}\big) - \frac{\epsilon}{2} (m_{c}+\theta_{c}).
\end{split}
\]
Using~\eqref{eq:1206} in~\eqref{eq:404} yields
\begin{align}
&\Delta(t) -V \sum_{c} \expect{h_{c}(\nu_{c}(t)) \mid H(t)} + V\sum_{nc} \theta_{c}\,\expect{\varphinc(t)\mid H(t)} \notag \\
&\leq B_{1} - V \sum_{c} h_{c}(r_{c}^{*}) + V \sum_{nc} \theta_{c}\,q_{n}^{(c)*} + \frac{V\epsilon}{2} \sum_{c} (m_{c}+\theta_{c}) \notag \\
&= B_{1} - V\sum_{c} g_{c}(r_{c}^{*}) + V\sum_{nc} \theta_{c}\,\lambdanc + \frac{V\epsilon}{2} \sum_{c} (m_{c}+\theta_{c}),  \label{eq:905}
\end{align}
where the last equality uses $h_{c}(x) = g_{c}(x) - \theta_{c}\, x$ and $r_{c}^{*} = \sum_{n} \big(\lambdanc-q_{n}^{(c)*}\big)$ in~\eqref{eq:213}. The next lemma is useful.

\begin{lem} \label{lem:403}
Define $\widetilde{\nu}_{c}(t) = \min[Z_{c}(t), \nu_{c}(t)]$ as the virtual data served in queue $Z_{c}(t)$ in slot $t$. If $Q \geq \numax$ then $\abs{\nu_{c}(t)-\tnu_{c}(t)} \leq \epsilon$ for all $t$ under the $\UORA$ policy.
\end{lem}
\begin{IEEEproof}[Proof of Lemma~\ref{lem:403}]
The value of $\nu_{c}(t)$ is the solution to the convex program~\eqref{eq:406}-\eqref{eq:407}. If $Z_{c}(t) \geq Q$, then we have $Z_{c}(t) \geq \numax \geq \nu_{c}(t)$ and $\tnu_{c}(t) = \nu_{c}(t)$. If $Z_{c}(t)<Q$, the problem~\eqref{eq:406}-\eqref{eq:407} reduces to
\begin{align}
\text{maximize} &\quad Vh_{c}(\nu_{c}(t)) - w\,\nu_{c}(t) e^{w(Q-Z_{c}(t))} \label{eq:408}\\
\text{subject to} &\quad 0\leq \nu_{c}(t)\leq \numax. \label{eq:409}
\end{align}
For all $x>\epsilon$, we have by mean-value theorem
\begin{align}
V \frac{h_{c}(x)-h_{c}(\epsilon)}{x-\epsilon} &= V h_{c}'(x_{\epsilon}), \quad \text{for some $x_{\epsilon} \in (\epsilon, x)$} \notag \\
&\leq 0 < w e^{w(Q-Z_{c}(t))}, \label{eq:1207}
\end{align}
where the first inequality results from the choice of $\theta_{c}$ that yields $h_{c}'(x) \leq 0$ for all $x\geq\epsilon$. Rearranging terms in~\eqref{eq:1207},  for all $x>\epsilon$ we have
\[
Vh_{c}(x) - w\,x e^{w(Q-Z_{c}(t))} < Vh_{c}(\epsilon) - w\,\epsilon e^{w(Q-Z_{c}(t))}.
\]
Therefore, the solution to the problem~\eqref{eq:408}-\eqref{eq:409} must satisfy $\nu_{c}(t) \leq \epsilon$. Since $0\leq \widetilde{\nu}_{c}(t) \leq \nu_{c}(t) \leq \epsilon$,  we conclude $\abs{\nu_{c}(t)-\tnu_{c}(t)} \leq \epsilon$.
\end{IEEEproof}

Now, from Lemma~\ref{lem:403} we have
\begin{equation} \label{eq:410}
h_{c}(\nu_{c}(t)) \leq h_{c}(\tnu_{c}(t)) + \epsilon (m_{c} + \theta_{c}).
\end{equation}
If $\nu_{c}(t)=\tnu_{c}(t)$ then this is trivial. Otherwise, if $\tnu_{c}(t) < \nu_{c}(t)$ then by mean-value theorem we obtain for some $y\in (\tnu_{c}(t), \nu_{c}(t))$ that
\[
\begin{split}
h_{c}(\nu_{c}(t)) &= h_{c}(\tnu_{c}(t)) + \big(\nu_{c}(t)-\tnu_{c}(t)\big) h_{c}'(y) \\
&\leq h_{c}(\tnu_{c}(t)) + \epsilon (m_{c} + \theta_{c}).
\end{split}
\]
Plugging~\eqref{eq:410} into~\eqref{eq:905} and rearranging terms yield
\begin{align}
&\Delta(t) -V \sum_{c} \expect{h_{c}(\tnu_{c}(t)) \mid H(t)} \notag \\
&\qquad+ V\sum_{nc} \theta_{c}\,\expect{\varphinc(t)\mid H(t)} \notag \\
&\leq B_{1} + \frac{3V \epsilon}{2} \sum_{c} (m_{c}+\theta_{c}) - V\sum_{c} g_{c}(r_{c}^{*}) + V\sum_{nc} \theta_{c}\,\lambdanc. \label{eq:411}
\end{align}
Define the time average
\[
\overline{\widetilde{\nu}_{c}(t)} \triangleq \frac{1}{t} \sum_{\tau=0}^{t-1} \expect{\widetilde{\nu}_{c}(\tau)}.
\]
Define $\overline{\varphinc(t)}$, $\overline{\widetilde{d}_{n}^{(c)}(t)}$, and $\overline{\widetilde{\mu}_{ad_{c}}^{(c)}(t)}$ similarly. Taking expectation and time average over $\tau\in\{0, \ldots, t-1\}$ in~\eqref{eq:411}, dividing by $V$, rearranging terms, and applying Jensen's inequality to the functions $h_{c}(\cdot)$, we get
\begin{multline} \label{eq:415}
\sum_{c} h_{c} \Big(\overline{\widetilde{\nu}_{c}(t)}\Big) + \sum_{nc} \theta_{c} \Big(\lambdanc - \overline{\varphinc(t)}\Big) \geq \\
\sum_{c} g_{c}(r_{c}^{*}) - \frac{B_{1}}{V} - \frac{\expect{L(H(0))}}{Vt} - \frac{3\epsilon}{2} \sum_{c}(m_{c}+\theta_{c}).
\end{multline}
Adding and subtracting $\sum_{c} \theta_{c}\, \overline{\widetilde{\nu}_{c}(t)}$ at the left-hand side of~\eqref{eq:415} and using the definition of $h_{c}(\cdot)$ yield
\begin{multline} \label{eq:417}
\sum_{c} g_{c} \Big(\overline{\widetilde{\nu}_{c}(t)}\Big) + \sum_{c} \theta_{c} \left\{\sum_{n} \lambdanc - \sum_{n} \overline{\varphinc(t)} -\overline{\widetilde{\nu}_{c}(t)} \right\} \\
\geq \sum_{c} g_{c}(r_{c}^{*}) - \frac{B_{1}}{V} - \frac{\expect{L(H(0))}}{Vt} - \frac{3\epsilon}{2} \sum_{c} (m_{c}+\theta_{c}).
\end{multline}
From~\eqref{eq:229} and~\eqref{eq:215} we have
\begin{gather*}
\Dnc(t+1) \geq \Dnc(t)-\varphinc(t) + \widetilde{d}_{n}^{(c)}(t), \\
Z_{c}(t+1) = Z_{c}(t) - \widetilde{\nu}_{c}(t) + \sum_{a} \widetilde{\mu}_{ad_{c}}^{(c)}(t).
\end{gather*}
Taking expectation and time average yields
\begin{gather} 
\overline{\varphinc(t)} \geq \overline{\widetilde{d}_{n}^{(c)}(t)} - \mathbb{E}\big[\Dnc(t)\big]/t, \label{eq:412} \\
\overline{\widetilde{\nu}_{c}(t)} = \sum_{a} \overline{\widetilde{\mu}_{ad_{c}}^{(c)}(t)} + \mathbb{E}\big[Z_{c}(0)-Z_{c}(t)\big]/t. \label{eq:413}
\end{gather}
Now, by the law of flow conservation, the sum of exogenous arrival rates must be equal to the sum of delivered throughput, time averages of dropped packets, and queue growth rates. In other words, we have for each class $c\in \cC$ and for all slots $t$
\begin{equation} \label{eq:414}
\sum_{n} \lambdanc = \sum_{n} \overline{\widetilde{d}_{n}^{(c)}(t)} + \sum_{a} \overline{\widetilde{\mu}_{ad_{c}}^{(c)}(t)} + \frac{1}{t}\sum_{n}  \mathbb{E}\big[\Qnc(t)\big].
\end{equation}
Combining~\eqref{eq:412}-\eqref{eq:414} yields
\begin{multline} \label{eq:416}
\sum_{n} \lambdanc - \sum_{n} \overline{\varphinc(t)} -\overline{\widetilde{\nu}_{c}(t)} \leq \\
\frac{1}{t} \expect{Z_{c}(t)-Z_{c}(0)} + \frac{1}{t} \sum_{n} \expect{\Dnc(t)+\Qnc(t)}.
\end{multline}

Finally, using the boundedness of queues $\Qnc(t)$, $\Dnc(t)$, and $Z_{c}(t)$ in Lemma~\ref{lem:502} and the continuity of $g_{c}(\cdot)$, we obtain from~\eqref{eq:413} and~\eqref{eq:416} that
\begin{gather}
\lim_{t\to\infty} g_{c}\Big(\overline{\widetilde{\nu}_{c}(t)}\Big) = g_{c}\Big(\lim_{t\to\infty} \sum_{a} \overline{\widetilde{\mu}_{ad_{c}}^{(c)}(t)}\Big) = g_{c}(\overline{r_{c}}), \label{eq:1210} \\
\lim_{t\to\infty} \bigg(\sum_{n} \lambdanc - \sum_{n} \overline{\varphinc(t)} -\overline{\widetilde{\nu}_{c}(t)} \bigg)\leq 0, \label{eq:1211}
\end{gather}
where the last equality of~\eqref{eq:1210} uses the definition in~\eqref{eq:504}. Taking a limit of~\eqref{eq:417} as $t\to\infty$ and using~\eqref{eq:1210} and~\eqref{eq:1211}, we obtain
\[
\sum_{c} g_{c} (\overline{r_{c}} ) \geq \sum_{c} g_{c}(r_{c}^{*}) - \frac{B_{1}}{V} - \frac{3\epsilon}{2} \sum_{c} (m_{c}+\theta_{c}),
\]
where the constant $B_{1}$, defined in~\eqref{eq:1209}, is
\begin{equation} \label{eq:1601}
\begin{split}
B_{1} &\triangleq \abs{\cN}\abs{\cC} \left[(\mumaxout + \dmax)^{2} + (\amax + \mumaxin)^{2} + 2\dmax^{2}\right] \\
&\quad + \abs{\cC} \Big[ w (2\deltamax + \epsilon)  +  e^{w(\numax + \mumaxin)} + \frac{w\epsilon}{2}e^{wQ} + e^{wQ}\Big].
\end{split}
\end{equation}
The proof is complete.

\section{} \label{appendix:203}

\begin{lem} \label{lem:101}
If a queue process $\{Q(t)\}$ satisfies
\begin{equation} \label{eq:117}
Q(t+1) \leq \big(Q(t) - b(t)\big)^{+} + a(t),
\end{equation}
where $a(t)$ and $b(t)$ are nonnegative bounded random variables with $a(t) \leq a_{\text{max}}$ and $b(t) \leq b_{\text{max}}$, then there exists a positive constant $B$ such that
\[
\frac{1}{2} \big(Q^{2}(t+1) - Q^{2}(t)\big) \leq B - Q(t)\,\big(b(t)-a(t)\big).
\]
\end{lem}
\begin{IEEEproof}[Proof of Lemma~\ref{lem:101}]
Squaring both sides of~\eqref{eq:117} yields
\begin{align*}
Q^{2}(t+1) &\leq \big(Q(t) - b(t)\big)^{2} + a^{2}(t) + 2\, Q(t)\,a(t) \\
&\leq Q^{2}(t) + B - 2\, Q(t)\, \big(b(t) - a(t) \big)
\end{align*}
where $B$ is a finite constant satisfying $B \geq a_{\text{max}}^{2} + b_{\text{max}}^{2}$. Dividing the above by two completes the proof.
\end{IEEEproof}

\end{document}